\documentclass[11pt,notitlepage,twoside]{article}
\usepackage{latexsym}
\usepackage{amsmath}
\usepackage{amsfonts}
\usepackage{amssymb}
\usepackage{amscd}
\usepackage{color}
\renewcommand{\a }{\alpha }
\renewcommand{\d}{\delta }
\newcommand{\D }{\Delta }

\newcommand{\e }{\varepsilon }
\newcommand{\g }{\gamma}

\renewcommand{\l }{\lambda }

\newcommand{\n }{\nabla }

\newcommand{\s }{\sigma }

\renewcommand{\o }{\omega }

\renewcommand{\O }{\Omega }

\newcommand{\ov}{\overline}
\newcommand{\intbar}{\mathop{\int\makebox(-13.5,0){\rule[4pt]{.7em}{0.3pt}}%
\kern-6pt}\nolimits}

\newcommand{\be}{\begin{equation}}
\newcommand{\ee}{\end{equation}}
\newcommand{\bes}{\begin{equation*}}
\newcommand{\ees}{\end{equation*}}
\newcommand{\ba}{\begin{eqnarray}}
\newcommand{\ea}{\end{eqnarray}}
\newcommand{\bas}{\begin{eqnarray*}}
\newcommand{\eas}{\end{eqnarray*}}
	
\newenvironment{pf}{\noindent{\sc Proof}.\enspace}{\rule{2mm}{2mm}\medskip}
\newenvironment{pfn}{\noindent{\sc \bf Proof }}{\rule{2mm}{2mm}\medskip}

\newcommand{\R}{\mathbb{R}}

\renewcommand{\o }{\omega }
\newcommand{\vn }{\varUpsilon }
\newcommand{\nto }{\nrightarrow}
\newcommand{\mS }{\mathcal{S}}

\author{ Mohamed Ben Ayed$^a$\thanks{  E-mail : {M.BenAyed@qu.edu.sa} } and
 Khalil El Mehdi$^{a,b}$\thanks{Corresponding author. E-mail : {K.Jiyid@qu.edu.sa} } \\ 
{\footnotesize
a : Department of Mathematics, College of Science, Qassim University, Buraydah 51542, Saudi Arabia}\\
{\footnotesize
 b : Facult\'e des Sciences et Techniques, Universit\'e de Nouakchott, Nouakchott, Mauritania.}\qquad\quad\quad\quad
}
\date{}

\title{\bf Multi-Bubble Blow-up Analysis  for an Almost Critical Problem} 
\begin{document}

\newtheorem{lem}{Lemma}[section]
\newtheorem{pro}[lem]{Proposition}
\newtheorem{thm}[lem]{Theorem}
\newtheorem{rem}[lem]{Remark}
\newtheorem{cor}[lem]{Corollary}
\newtheorem{df}[lem]{Definition}

\maketitle

\noindent{\bf Abstract:} Consider a smooth, bounded domain $\O\subset \mathbb{R}^n$ with $n\geq 4$ and a smooth positive function 
$V$. We analyze the asymptotic behavior of a sequence of positive solutions $u_\e$ to the equation $-\Delta u +V(x)u =u^{\frac{n+2}{n-2}-\e}$ in $\O$ with zero Dirichlet boundary conditions, as $\e\to 0$. We determine the precise blow-up rate and characterize the locations of interior concentration points in the general case of multiple blow-up, providing an exhaustive description of interior blow-up phenomena of this equation. Our result is established through a delicate analysis of the gradient of the corresponding Euler-Lagrange functional.

\bigskip


\bigskip

\noindent {\bf MSC $2020$}: 35A15, 35J20,  35J25.

%%%%%%%%%%%%%%%%%%%%%%%%%%%%%%%%
\section{Introduction and Main Results}
%%%%%%%%%%%%%%%%%%%%%%%%%%%%%%%%
In this paper, we are interested in the behavior of solutions to the following nonlinear  problem
\begin{align}
(\mathcal{P}_{\e, V}):\qquad 
\begin{cases}
-\Delta u +V u = u^{p-\varepsilon}\quad &\mbox{in}\quad \Omega,
\\
\quad  u > 0 \quad &\mbox{in}\quad \Omega,
\\
\quad u =0\,\quad &\mbox{on}\quad \partial\Omega,
\end{cases}
\end{align}
where $\Omega$ is a smooth bounded domain in $\mathbb{R}^n$, $n\geq 4$, the potential $V$ is a $C^3$ positive function on $\overline{\Omega}$, $p+1= (2n)/(n-2)$ is the critical Sobolev exponent for the embedding $H^1_0(\Omega) \hookrightarrow L^q(\Omega)$ and $\varepsilon$ is a small positive parameter.  

The problem in the form of $(\mathcal{P}_{\e, V})$ appears in various physical models, including non-relativistic Newtonian gravity and quantum transport, as discussed in \cite{BGDB, BLJS, Pe} and their respective references. In addition, the understanding of the behavior of the solutions to this problem plays a crucial role in issues related to prescribed scalar curvature problems in differential geometry, see for instance \cite{DH} and the references therein.

Over the past decades, a substantial body of literature has developed on the blow-up phenomena of solutions to semilinear equations with critical exponents, a topic too broad to summarize fully. Here, we highlight only a selection of some works, from which we hope a more comprehensive bibliography can be constructed. In the case where $V\equiv 0$, the asymptotic behavior of solutions that blow up at a single point has been well-established in \cite{BP, R1, Ha}. The scenario of multiple concentration points, still for $V\equiv 0$, was later investigated in \cite{R2, BLR, R3}. For $V\not\equiv 0$, single-blow-up behavior in the interior has been studied in the three-dimensional case by \cite{FKK}. The focus of this paper is on multiple blow-up points in the interior for dimensions greater than or equal to 4. Other blow-up scenarios have also been studied, including those for sign-changing solutions and for  positive solutions under Neumann boundary conditions. For additional references, we refer the reader to \cite{CP, KL1, KL, DHR, H, MM, AB1, ABE, yyli, BM, KMS} and the works cited therein. For the inverse problem of constructing blow-up solutions in this setting, we direct the reader to \cite{2ABE, Clus, BA, BLR, EM, BEM, MRV, RW1, RW2, DDM, MS, Pi} and their respective references.
 
%%%%
Notice that problem $(\mathcal{P}_{\e, V})$ has a variational structure, with its solutions being the positive critical points of the functional
\be\label{eq:1}
I_{\e, V}(u)=\frac{1}{2}\int_{\O}\left(|\n u|^2+V(x)u^2\right)-\frac{1}{p+1-\e}\int_{\O}|u|^{p+1-\e} 
\ee
defined on  $H^1_0(\O)$  equipped with the norm $\|.\|$ and its corresponding inner product given by
\be\label{eq:4}
\left\langle v,w\right\rangle:=\int_{\O}\n v\cdot\n w+\int_{\O}Vvw \quad ; \quad \|w\|^2 :=\int_{\O}|\n w|^2+\int_{\O}V w^2.
\ee
Since $V$ is a continuous positive function on $\ov{\O}$, we see that this norm is equivalent to the two norms $||.||_{0}$ and $||.||_{1}$ of $H^1_0(\O)$ and $H^1(\O)$ respectively.\\
Note that all solutions $u_\e$ of $(\mathcal{P}_{\e, V})$ satisfy $\|u_\e\| \geq C$ with a positive constant $C$ independent of the parameter $\e$. Therefore, the concentration compactness principle (\cite{S, L}) implies that if $u_\e$ is an energy-bounded solution of $(\mathcal{P}_{\e, V})$ that weakly ( but not strongly) converges to  $ \o $ (where $ \o = 0 $ or $ \o $ is a solution of $(\mathcal{P}_{0, V})$), then $u_\e$  must concentrate at a finite number $N$ of points in $\overline{\O}$. More precisely, $u_\e$ can be expressed as:
\begin{align}
 & u_{\e} = \o + \sum_{i=1}^N \d_{{a_{i,\e}}, \l_{i,\e}} +v_\e  \quad \mbox{ where } \label{eq:3} \\
& \left\|v_{\e}\right\|  \rightarrow 0,\, \l_{i, \e}\,d(a_{i,\e},\partial\O) \rightarrow \infty,\,a_{i,\e} \to \overline{a}_i\in \overline{\O} \quad\forall i \quad \text { as } \quad \e \rightarrow 0,\label{eq:bis}\\
&  \e_{i j}:=\left(\frac{\l_{i, \e}}{\l_{j, \e}}+\frac{\l_{j, \e}}{\l_{i, \e}}+\l_{i, \e} \l_{j, \e}\left|a_{i, \e}-a_{j, \e}\right|^2\right)^{\frac{2-n}{2}} \rightarrow 0 \text { as } \e \rightarrow 0 \quad\forall i \neq j \label{eq:bbis}
\end{align}
and where $\d_{a,\l}$ are the functions (called bubbles) defined by 
\be\label{bubble}
   \d_{a,\l } (x) = c_0\frac{\l^{\frac{n-2}{2}}}{(1+\l^2\mid
 x-a\mid^2)^{\frac{n-2}{2}}}, \qquad \l>0, a,\,x \in \R^n,\, c_0=(n(n-2))^{\frac{n-2}{4}} 
 \ee
and which are, see \cite{CGS},  the only solutions of the following problem 
$$
 -\D u= u^\frac{n+2}{n-2} ,\, u > 0\,\quad \mbox{ in } \quad \R^n.
$$

 It is worth noting that the following related problem was investigated in \cite{KL1} for the case $ n \geq 4$ 
\begin{align*}\label{8}
(\mathcal{Q}_{\e, V}):\qquad \begin{cases}
-\Delta u +\e V u = u^{p} , \quad   u > 0 \quad &\mbox{in}\quad \Omega, \\
\quad u =0\,\quad &\mbox{on}\quad \partial\Omega . 
\end{cases}
\end{align*}
In \cite{KL1}, by considering a blowing-up family of solutions $ ( u_\e ) $, the authors proved that when $ V < 0 $, the concentration points $a_{i,\e}$'s remain away from both the boundary and each other, that is, $ d(a_{i, \e} , \partial \O ) \geq c > 0 $ and $ | a_{i, \e} - a_{j, \e} | \geq c > 0$ for each $ i \neq j $. In contrast, for our problem, $ ( \mathcal{P}_{\e , V } ) $, the behavior is fundamentally different. We believe that when $ V < 0 $ and $-\D +V$ is coercive, blow-up phenomena do not occur, whereas they may arise when $ V > 0 $. Moreover, some concentration points may converge to the same location (see the theorems below). Additionally, it appears possible to construct solutions featuring concentration points -either isolated or clustered- located on the boundary, see the end of Subsection 2.1. We plan to explore this question in future work.

In this paper, our goal is to examine the asymptotic behavior of interior concentrating solutions $u_\e$, as $\e\to 0$, for problem $(\mathcal{P}_{\e, V})$. More specifically, we focus on the case where
\be\label{int}
d(a_{i,\e}, \partial\O) \geq d_0>0 \quad \forall\,i\in\{1,\cdots, N \}.
\ee

Notice that, the existence of solutions to $(\mathcal{P}_{\e, V})$ can be proved by noting that the infimum
$$
\inf \,\{\int_\O \left(|\n u |^2 +V\,u^2\right):\quad u\in H^1_0(\O) \,  \mbox{ and } \,  \int_\O |u|^{p+1-\e}=1\, \}
$$
is achieved, owing to the compactness of the embedding $H^1_0(\Omega) \hookrightarrow L^q(\Omega)$ for $ 1 \leq q <  2n/(n-2)$. 
 
 %%%%%%%%%%%%%%%%%%%%%%%%%%%%%%%%

%%%%%%%%%%%%%%%%%%%%%%%%%%%%%%%%
The focus of this paper is to more precisely describe the blow-up behavior of solutions to $(\mathcal{P}_{\e, V})$. Our goal is to determine the exact blow-up rate and characterize the locations of interior blow-up points, offering a comprehensive description of the interior blow-up phenomena for this problem. More specifically, we prove:
\begin{thm}\label{th:t1}
 Let $n\geq 4 $ and  $0 < V \in C^3(\ov{\O})$ be a positive function and let $(u_\e)$ be a sequence of solutions of $(\mathcal{P}_{\e, V})$ having the form \eqref{eq:3}, with $ \o = 0$ (i.e. $ u_ \e $ converges weakly to $ 0 $), and satisfying \eqref{eq:bis}, \eqref{eq:bbis} and \eqref{int}. Then, the following results hold:
\begin{enumerate}
  \item[(a)]
  Each  concentration point  $ a_{i,\e} $  converges to a critical point  $z_i$ of $ V $. In addition,  there exists a dimensional constant $\kappa_1(n)  > 0$  (see \eqref{RES1} for the precise value) such that
 $$
 \l_{i,\e} ^2 \, = \, \kappa_1 (n) \, V(z_i) \, \frac{ ( | \ln  \e | / 2 )^{ \s_n}  }{\e } (1+o_\e(1)), 
 $$
 where $ \s_4 := 1 $ and $ \s_n := 0 $ for $ n \geq 5 $.
  \item[(b)] For  a non degenerate critical point $z$ of $ V$, let $ B(z):= \{ a_{i, \e}: \lim | a_{i, \e} - z | = 0 \}$.   There are three alternatives
  \begin{enumerate}
    \item[(i)] Either,  $ B(z) = \emptyset$, in this case $z$ is not a blow up point of the family $(u_\e)$,  
        \item[(ii)] or $ \# B(z) = 1 $, let $ B(z)= \{ a_{i, \e} \}$, it holds that $ \beta ( \e ) | a_{i, \e} - z | $ is uniformly bounded where 
        $$ \beta ( \e ) := ( | \ln \e | \, \mbox{ if } n = 4 \, ; \, \, \,  \e^{-1/2} | \ln \e |^{-1} \, \mbox{ if } n = 5 \, ; \, \, \,  \e^{-1/2}  \, \mbox{ if } n \geq 6) . $$ In this case, $z$ is an  \emph{isolated simple blow up point} for the family $(u_{\e})$.
    \item[(iii)]  or  $ \# B(z) \geq 2 $, let $ B(z)= \{ a_{i_1, \e},\cdots,a_{i_m, \e} \}$, it holds that $  | \ln \e | ^{ \s_n / 4 } \e^{-(n-4)/(2n) }  | a_{i_j,\e} -  z  | $  is uniformly bounded for each $ j \in \{ 1, \cdots, m \}$. Furthermore, there exist at least $ m-1$ points $ a_{i_j, \e}$'s such that $  | \ln \e | ^{ \s_n / 4 } \e^{-(n-4)/(2n) }  | a_{i_j,\e} -  z | $ is uniformly bounded below. 
    In this case $z$ is a  \emph{non simple blow up point} for the family $(u_{\e})$. 
  \end{enumerate} 
  \item[(c)] Let $z  \in \O$ be  a non degenerate critical point of $ V $ which is a non simple  blow up point  and  let $ B(z) = \{a_{1,\e}, \cdots,a_{m,\e}\}$.        Denoting
      \begin{equation}\label{bi} 
        b_{i,\e}:= \ \kappa_2(n) \, V (z)^{(n-4)/(n-2)} \, \eta (\e) ^{-1} \, ( a_{i,\e} - z )   \quad  \mbox{ for } 1\leq i\leq m,
      \end{equation}
  where 
 \be \label{etaeps} \eta(\e):= (  \e^{(n-4)/(2n)} \, \,  \mbox{ if } n \geq 5 \, \, ; \, \, \,  \left( 2 / | \ln \e | \right)^{1/4} \, \,  \mbox{ if } n = 4 ) \ee     
and  $\kappa_2(n) > 0$ is a dimensional constant (see \eqref{kappa2} for the precise value),  we have that
       $(b_{1,\e},\cdots,b_{m,\e})$ converges to $(\ov{b}_1,\cdots,\ov{b}_m)$ which is a critical point of the following Kirchhoff-Routh type  function:
\begin{align} \label{FYN} & \mathcal{F}_{z,m}:  \mathbb{F}_m( \O )\to \R ; \quad    \mathcal{F}_{z,m}(\xi_1,\cdots,\xi_m):= \sum_{i=1}^{m} D^2 V (z)( \xi_i , \xi_i )  -  \sum_{1 \leq i ,  j \leq m}\frac{1}{| \xi_i - \xi_j | ^{n-2}},  \end{align}
where $ \mathbb{F}_m(  \O ):=  \{   (\xi_1,\cdots,\xi_m) \in \O^m : \, \,  \, \xi_i \neq \xi_j   \,  \mbox{ for } i \neq j  \} $.
\end{enumerate}
\end{thm}

The next theorem provides a description of the interior blow-up phenomenon occurring when the sequence of solutions converges weakly to a non-zero limit.
\begin{thm}\label{th:t1om1}
 Let $n\geq 7$ and  $ 0 < V \in C^3(\ov{\O})$ be a positive function and let $(u_\e)$ be a sequence of solutions of $(\mathcal{P}_{\e, V})$ having the form \eqref{eq:3} with $ \o \neq 0$ and satisfying \eqref{eq:bis}, \eqref{eq:bbis} and \eqref{int}. Then, the conclusions of Theorem \ref{th:t1} hold with, in case $(ii)$, $ \beta(\e) = -1/2 $ for $ n \geq 8 $ and $ \beta(\e) = -1/4 $ for $ n = 7 $.
\end{thm}
\begin{rem}
In lower dimensions, that is, for $3\leq n\leq 5$, we believe that $(\mathcal{P}_{\e, V})$ admits no interior blowing up solutions whose weak limit is a non-zero solution of the limit problem $(\mathcal{P}_{0, V})$, see $(2)$ in Remark \ref{rem om2}.   This question will be addressed in a future work.
\end{rem}

The proofs of our results are based on some balancing conditions satisfied by the concentration parameters, that is,  relationships that ensure the equilibrium between different factors governing the blow-up behavior of the solutions. These conditions are derived by performing an asymptotic expansion of the gradient of the functional $I_{\e, V}$ and testing the equation with suitable vector fields, which leads to constraints on the concentration points and the corresponding rates at which the solution blows up. Through a careful study of these conditions, we obtain a precise description of the blow-up phenomenon, including the locations and rates of concentration.

The remainder of the paper is organized as follows: In Section 2, we present the balancing conditions that relate the concentration parameters. Section 3 is devoted to proving the estimates for the concentration rates. Finally, in Section 4, we establish the proof of our result.

%%%%%%%%%%%%%%%%%%%%%%%
\section{ Balancing conditions of the concentration parameters } 
%%%%%%%%%%%%%%%%%%%%%%%%
In this section, we aim to identify the relationships that maintain equilibrium between the various factors governing the blow-up behavior of the solutions of $(\mathcal{P}_{\e, V})$. To achieve this goal, let $(u_\e)$ be a sequence of solutions of $(\mathcal{P}_{\e, V})$ having the form \eqref{eq:3} and satisfying \eqref{eq:bis}, \eqref{eq:bbis} and \eqref{int}. 
We notice that the function $ \d_{a,\l} $, introduced in \eqref{bubble}, does not belong to $ H_0^1( \O ) $, thus we need to introduce its projection onto this space. Different projections are introducing in some works. In our case, we will use the following projection: for $a\in\O$ and $\l>0$, 
\begin{equation} \label{x}
\begin{cases}
(-\Delta +V) (\pi \delta_{a,\lambda}) = \delta_{a,\lambda}^{ (n+2) / (n-2) }\quad &\mbox{in}\quad \Omega,  \\
\quad\pi \delta_{a,\lambda} =0\,\quad &\mbox{on}\quad \partial\Omega . 
\end{cases}
\end{equation}
Observe that, following \cite{BCtopology}, we can modify the variables $ a_{i, \e}$'s and $ \l_{i, \e}$'s (using the same notation for the new variables) so that the decomposition of $u_\e$ becomes as follows
\begin{align} 
 &  u_\e =  \a_{0, \e} \o + \sum_{i=1 }^N \a_{i, \e} \pi \d_{a_{i,\e},\l_{i,\e}} + v_\e,  \qquad \mbox{ with } \label{eq:6} \\  
  & \a_{i, \e} \to 1 , \, \, d(a_{i, \e} , \partial \O ) \geq 2 d_0 , \, \,   \l_{i,\e}  \to \infty,  \, \, \e_{ij } \to 0 , \, \, \, \| v_\e \| \to 0 \, \,  \mbox{ and } \, \,  v _\e \in E_{a, \l}^\perp \quad  \mbox{ where } \label{assump2} \\ 
   & E_{a,\l} := \mbox{span} \Big\{ \o,  \pi \d_{a_i, \l_i} , \frac{ \partial (\pi \d_{a_i, \l_i}) }{ \partial \l_i } ,  \frac{ \partial (\pi \d_{a_i, \l_i}) }{ \partial a_{i,j}} , \, \, i \in \{1, \cdots , N \} , \, \, j \in \{1, \cdots , n \} \Big\}.  \label{Eal} 
\end{align}
We precise that the orthogonality is taken with respect to the inner product defined in \eqref{eq:4}. 

We note that, for the sake of simplicity in presentation, we will omit the index $ \e $ from the variables $ a_{i,\e} $, $ \l_{i,\e} $ and $ \a_{i,\e} $ in some of the proofs that follow. We begin this section with the following result, which will be crucial for the computations that follow.
\begin{pro} \label{epslambda}
For each $i \in \{1, \cdots, N \}$, it holds that: $ \e \ln \l_{i, \e} \to 0 . $
\end{pro}
\begin{pf} Multiplying $ ( \mathcal{P}_{\e, V} ) $ by $ \pi \d_{a_i, \l_i} $ and integrating over $ \O$, we obtain 
$$ \a_i \| \pi \d_{a_i, \l_i} \|^2 + \sum_{j \neq i } \a_j \langle \pi \d_{a_j, \l_j}, \pi \d_{a_i, \l_i} \rangle + \a_0 \langle \o , \pi \d_{a_i, \l_i} \rangle  = \int_\O u_\e ^{ p-\e} \pi \d_{a_i, \l_i} . $$
 Observe that 
\be \label{delta w } \langle \o , \pi \d_{a_i, \l_i} \rangle  = \int_{ \O } \d_{a_i, \l_i} ^{ (n+2)/(n-2)} \o = \ov{c}_1 \frac{ \o(a_i)}{ \l_i ^{(n-2) / 2 } } + O \Big( \frac{ \ln \l_i  } { \l_i ^{ (n+2) / 2 } }  \Big) \ee
 with $$  \ov{c}_1 := c_0^{\frac{n+2}{n-2}} \int_{\R^n } \frac{1}{ ( 1+ | x |^2 )^{(n+2)/2} }dx. $$

From equations $(31)$ and $(5)$ of \cite{Clus}, we have 
\be \label{aqw1} \a_i \| \pi \d_{a_i, \l_i} \|^2 = \a_i S_n + o(1) \qquad \qquad  \mbox{ and } \qquad\qquad  0 \leq  \pi \d _{ a_i, \l_i } \leq \d _{ a_i, \l_i } . \ee
In addition, for $ j \neq i$, using Estimate 1 of \cite{Bahri-book}, it follows that 
\be \label{1wxc}
 \langle \pi \d_{a_j, \l_j}, \pi \d_{a_i, \l_i} \rangle = \int_\O \d_{a_j, \l_j}^{(n+2)/(n-2)} \pi \d_{a_i, \l_i} \leq  \int_\O \d_{a_j, \l_j}^{(n+2)/(n-2)} \d_{a_i, \l_i} \leq c \e_{ij} = o(1) . \ee
Furthermore, using the standard analysis formulae: for $ \g > 0 $, $ s , t , t_i \in \R $, 
$$ | s+ t |^\g = | s |^\g + O( | s |^{\g - 1 } | t | +  | t | ^{\g  } ) \qquad \qquad  \mbox{ and } \qquad \qquad \Big| \sum t_i \Big|^\g \leq c \sum | t_i |^\g , $$
we deduce that 
\begin{align*} \int_\O u_\e ^{ p-\e} \pi \d_{a_i, \l_i} = \a_i ^{p-\e} \int_\O   (\pi \d_{a_i, \l_i} )^{ p-\e + 1} + O \Big(  \int_\O   (\pi \d_{a_i, \l_i} )^{ p-\e} \Big[  \o + \sum_{j \neq i }  (\pi \d_{a_j, \l_j} ) + | v _\e | \Big] \Big)  \\
+ O \Big( \int_\O   (\pi \d_{a_i, \l_i} ) \Big[  \o ^p + \sum_{j \neq i }  (\pi \d_{a_j, \l_j} )^{ p-\e} + | v _\e |^{ p-\e} \Big] \Big) . 
\end{align*}
Now, using \eqref{aqw1} and the fact that: for $ 1 \leq \g \leq p $ and $ t > 0$, it holds: $ t^{\g } \leq t + t^p $, we deduce that 
\begin{align} 
 &  \int_\O   (\pi \d_{a_\ell, \l_\ell} )^{ p-\e}   (\pi \d_{a_k, \l_k} ) \leq \int_\O   \d_{a_\ell, \l_\ell} ^{ p-\e}   \d_{a_k, \l_k}  \leq  \int_\O  [ \d_{a_\ell, \l_\ell} + \d_{a_\ell, \l_\ell}^{ p} ]  \d_{a_k, \l_k}  = o(1) , \notag \\
 &  \int_\O   (\pi \d_{a_i, \l_i} )^{ p-\e}  \o +  \int_\O   (\pi \d_{a_i, \l_i} )  \o ^p    \leq  c  \int_\O  [ \d_{a_i, \l_i} + \d_{a_i, \l_i}^{ p} ]   = O \Big(\frac{1}{ \l_i ^{(n-2)/2}}\Big)  = o(1),  \label{om9}\\
 & \int_\O   (\pi \d_{a_i, \l_i} )^{ p-\e}  | v _\e | + \int_\O   (\pi \d_{a_i, \l_i} )  | v _\e |^{ p-\e} \leq \int_\O  [ \d_{a_i, \l_i}  + \d_{a_i, \l_i} ^{ p} ]  | v _\e | + \int_\O   \d_{a_i, \l_i}  [ | v _\e | +   | v _\e |^{ p} ] = o(1) ,  \notag
\end{align}
by using Holder's inequality and the fact that $ \| v_\e \| \to 0 $ for the last equation and \eqref{1wxc} for the first one. Finally, using equations $(64)$ of \cite{EM} and $(30)$ of \cite{Clus}, we deduce that 
$$  \int_\O   (\pi \d_{a_i, \l_i} )^{ p-\e + 1} =  \int_\O    \d_{a_i, \l_i} ^{ p-\e + 1} + O \Big(  \int_\O   \d_{a_i, \l_i} ^{ p-\e } | \d_{a_i, \l_i} - \pi \d_{a_i, \l_i} | \Big) = c_0^{-\e} \l_{i} ^{-\e (n-2) / 2 } S_n  + o(1) . $$ 
Combining the previous estimates, we derive the result. 
 \end{pf}

 \subsection{ Case of zero weak limit }

In this subsection, we will assume that $ \o = 0 $. 

Next, we will present the balancing conditions satisfied by the variables  $\a_{i, \e}$'s, $ \l_{i, \e}$'s and $ a_{i, \e}$'s, as defined in equation \eqref{eq:6}. To this end, note that  \eqref{assump2} and Proposition \ref{epslambda} imply that 
$$ (\a_\e, a_\e , \l_\e) :=  ( \a_{i, 1}, \cdots, \a_{i, N}, a_{i, 1}, \cdots, a_{i, N},  \l_{i, 1}, \cdots, \l_{i, N} ) \in \mathcal{O} ( N, \mu) $$ 
for small positive real $ \mu$, where the set $\mathcal{O} ( N, \mu)$ is defined by 
 \begin{align}
 \mathcal{O} ( N, \mu) := \{ (\a, a, \l) \in (0,\infty)^N \times \O ^N \times ( \mu^{-1} , \infty)^N : \, \, & | \a_i - 1 | < \mu ; \, \, d(a_i, \partial \O ) > c ;  \label{eq:2} \\
  &  \e \ln \l_i < \mu ; \, \, \forall \, \, i  \, \mbox{ and }\, \,  \e_{ij} < \mu  \, \, \forall \, \, i \neq j \}  \nonumber 
  \end{align}
with $ c $ is a fixed positive constant. \\
This implies that we are in a position to apply all the asymptotic estimates established in \cite{Clus}. Thus, letting  $ \underline{u}_\e := \sum_{i=1}^N \a_{i, \e} \pi \d_{a_{i,\e}, \l_{i, \e}} $, we see that the function $ v_\e $ satisfies equation $ (27) $ of \cite{Clus} which has a unique solution. Therefore, applying Proposition 2.2 of \cite{Clus}, we derive that  $ v_\e $  satisfies the following estimate 
\be \label{estveps} 
\| v _\e \| \leq c \e + c \, \sum \frac{ \ln^{ \s_n} (\l_i) }{ \l_i ^2 } + c \,  \sum_{ j \neq i } \begin{cases}
\e_{ij} & \mbox{ if } n \leq 5 , \\
\e_{ij}^{(n+2) /( 2(n-2)) } \ln( \e_{ij}^{-1}) ^{(n+2)/(2n)} & \mbox{ if } n \geq 6.
\end{cases} 
\ee 
In addition, the variables $\a_{i, \e}$'s, $ \l_{i, \e}$'s and $ a_{i, \e}$'s satisfy the following balancing conditions 
\begin{pro} \label{balancing}
For $ \e $ small, the following facts holds 
\begin{align}
& ( EAl_i) \quad \left( 1-\alpha_{i}^{p-1-\varepsilon}\lambda_{i}^{ -\varepsilon (n-2) / 2 }\right) = O(R_{\alpha_{i}}) \qquad \mbox{ with }   R_{\alpha_{i}} :=  \varepsilon + \frac{ \ln^{ \s_n} (\l_i) }{ \l_i ^2 } +\left\| v\right\|^{2}+\sum_{i\neq j}\varepsilon_{ij} , \label{EAli1} \\
& 	( EL_i)  \quad c_2    \e - \ov{ c }_2 \sum_{ j \neq i }  \l_{i,\e} \frac{ \partial \e_{ij} }{ \partial \l_{i, \e} } - c(n) \frac{ \ln ^{\s_n} ( \l_{i, \e} ) }{ \l_{i, \e} ^2 } V(a_{i, \e} ) \notag \\
&  \qquad = o\Big( \e + \sum_{ j \neq i } \e_{ij}  +\frac{ \ln ^{\s_n} ( \l_{i, \e} ) }{ \l_{i, \e} ^2 } \Big) + O \Big( \| v_\e\|^2 + \sum_{k \neq j } \e_{kj}^{\frac{n}{n-2}} \ln ( \e_{kj}^{-1}) + \sum_{j \neq i } \frac{1}{ ( \l_j\l_i) ^{\frac{n-2}{2}}} \Big) , \label{ELi1} \\   
&(EA_i)  \quad   c_2(n) \a_{i, \e} \frac{ \ln ^{\s_n} ( \l_{i, \e} ) }{ \l_{i, \e} ^3 } \n V(a _{i, \e} ) -  \ov{ c }_2 \sum_{ j \neq i } \a_{j,\e} \frac{1}{\l_{i,\e} }    \frac{ \partial \e_{ij} }{ \partial a_{i, \e} } = O ( R' _{a_i} ), \label{EAi1}  
  \end{align}
where  	
\begin{align*}
 &  c_2(4)   := \frac{1}{2} c_{0}^{2}\, meas(\mathbb{S}^{3}) \qquad ,    \qquad c_2(n):=\frac{n-2}{n}c_{0}^{2}\int_{\mathbb{R}^{n}}\frac{|x|^{2}}{(1+|x|^{2})^{n-1}} \quad\mbox{if n}\geq 5, \\
 & c(4)  :=  c_{0}^{2}\,mes(\mathbb{S}^{3}) \quad , \qquad c(n):=\frac{n-2}{2}c_{0}^{2}\int_{\mathbb{R}^{n}}\frac{|x|^{2}-1}{(1+|x|^{2})^{n-1}}> 0 \quad\mbox{if n}\geq 5, \\
 &  c_2  := \Big(\frac{n-2}{2} \Big) ^{2} c_{0}^{\frac{2n}  { n-2 }} \int_{\mathbb{R}^{n}} \frac{(\left|x\right|^{2}-1)}{(1+\left|x\right|^{2})^{n+1}} \ln(1+\left|x\right|^{2})dx  > 0 \,  ,  \quad \ov{c}_2   := c_0^{\frac{2n}{n-2}} \int_{\R^n } \frac{1}{ ( 1+ | x |^2 )^{(n+2)/2} }dx,\\
 &R'_{a_i} =   R_{ a_{i, \e} } + R_{ \a_{i, \e} } \Big[  \frac{ \ln ^{\s_n} ( \l_{i, \e} ) }{ \l_{i, \e} ^3 } + \sum _{j \neq i }  \frac{1}{\l_{i,\e} }  \Big|  \frac{ \partial \e_{ij} }{ \partial a_{i, \e} } \Big| \Big],\\
 & R_{a_i} =   \sum_{ j \neq i}  \l_j | a_i - a_j | \e_{ ij }^{\frac{n+1}{n-2}} + \left\| v\right\|^{2} + \frac{1}{\lambda_{i}^{4}} + (\mbox{if}\; n=5)\frac{ \ln\lambda_{i}}{\lambda_{i}^{4}}  \\
   &\quad  +\sum  \frac{1}{\lambda_{k}^{n-1}}  + \sum_{k \neq j } \e_{kj}^{n/(n-2} \ln \e_{kj}^{-1}  + \sum  \e_{ij}  \Big( \e + \frac{1}{\l_i} + \frac{1}{\lambda_{j}^{3/2}} \Big) + \e\frac{\ln^{\sigma_{n}} \l_{i}}{\l_{i}^{2}}.
  \end{align*}
  \end{pro}
  \begin{pf} Using the fact that $ u_\e $ is a critical point of $ I_{\e, V}$, we see that Propositions 3.3, 3.4 and 3.5 of \cite{Clus} imply that 
  \begin{align}
& \left( 1-\alpha_{i}^{p-1-\varepsilon}\lambda_{i}^{ -\varepsilon (n-2) / 2 }\right) = O(R_{\alpha_{i}}),\label{en1}\\
& 	\sum_{j\neq i} \a_j\, \ov{c}_2 \l_i \frac{\partial \e_{ij}}{\partial \l_i} \Big(1- \a_i^{p-1-\e} \, \l_i^{ -\e (n-2) / {2}} - \a_j^{p-1-\e} \, \l_j^{ - \e (n-2) / {2}} \Big) \notag\\
& \qquad \qquad  +  c_{0}^{-\varepsilon}\lambda_{i}^{-\varepsilon \frac{n-2}{2}}c_{2}\alpha_{i}^{p-\varepsilon}\varepsilon   
   - c(n) \alpha_{i}\frac{\ln^{\sigma_{n}}(\lambda_{i})}{\lambda_{i}^{2}} V(a_i) \left ( 2\alpha_{i}^{p-\varepsilon -1} c_{0}^{-\varepsilon }\lambda_{i}^{-\varepsilon \frac{n-2}{2}}-1 \right )  =  O ( R_{\l_i} ),\label{en2} \\
   &     c_2(n) \alpha_{i}\frac{\ln^{\sigma_{n}}(\lambda_{i})}{\lambda_{i}^{3}} \n V(a_i) \left ( 2\alpha_{i}^{p-\varepsilon -1} c_{0}^{-\varepsilon }\lambda_{i}^{-\varepsilon \frac{n-2}{2}}-1 \right )\notag \\
 &   \qquad \qquad  + \ov{c}_2 \sum_{ j \neq i } \a_j \frac{1}{ \l_{i} }  \frac{\partial \e_{ ij }}{\partial a_{i}} \Big( 1 - c_0^{-\e}  \l_i ^{ -\e \frac{n-2}{2} }  \a_i^{p-\e-1 } - c_0^{-\e}  \l_i ^{ -\e \frac{n-2}{2} } \a_j^{p-\e-1 }    \Big)  =  O ( R_{a_i}),\label{en3}
  \end{align}
  with 
 \begin{align*}
    R_{\l_i}  = &  \varepsilon^{2}+\left\| v\right\|^{2}  + \frac{1}{\lambda_{i}^{n-2}}+\frac{1}{\lambda_{i}^{4}} + (\mbox{if}\; n=6)\frac{ \ln\lambda_{i}}{\lambda_{i}^{4}}   +  \sum_{k \neq j }\e_{kj}^{\frac{n}{n-2}} \ln (\e_{kj}^{-1}) \\
 & \qquad  \qquad \qquad +  \sum_{j \neq i } \Big( \frac{ 1} { (\l_j  \l_i ) ^{ (n-2) / 2 } } +\e_{ij} \Big(  | a_i - a_j | ^2 | \ln | a_i - a_j | | ^{\s_n}  +  \frac{ 1} { \l_j ^{3/2}} +  \frac{ 1} { \l_i ^{3/2}}  \Big) \Big). 
  \end{align*}
 Putting \eqref{en1} in \eqref{en2} and \eqref{en3}, we get the desired estimates.
 \end{pf}

Observe that Equation $ (EL_i)$ in Proposition \ref{balancing} implies that, when $n\geq 4$, $V<0$, and the concentration points are away from the boundary, the following holds:
$$
 c'_1 \e  + c'_2 \sum_{ j \neq i }  \e_{ij} +  c'_3 \sum_{ i = 1 }^N \frac{ \ln ^{ \s_n} ( \l_{i, \e} )}{ \l_{i, \e} ^2} = o \Big( \e  +  \sum_{ j \neq i }  \e_{ij} + \sum_{ i = 1 }^N \frac{ \ln ^{ \s_n} ( \l_{i, \e} )}{ \l_{i, \e} ^2} \Big). 
 $$ 
 This leads to a contradiction, and thus we obtain the following remark, which highlights a fundamental difference between problems $(\mathcal{P}_{\e, V})$   and $(\mathcal{Q}_{\e, V})$.
 \begin{rem}\label{rn}
 When $n\geq 4$, $V<0$ and $-\D +V$ is coercive, problem $(\mathcal{P}_{\e, V})$ admits no solutions that blow up solely at interior points as $\e$ goes to zero.
 \end{rem} 
 Before concluding this subsection, we briefly comment on Assumption \eqref{int}, which states that all concentration points lie in the interior of the domain. 
To gain a better understanding of the phenomenon, we consider the case   $ N = 1 $.
 In this setting, when the concentration point $ a_\e $  is close to the boundary, additional terms appear in the equations of Proposition \ref{balancing}. It can be shown that the resulting equations in this case are given by, for $ n \geq 4 $, 
 \begin{align}
& ( EAl') \quad \left( 1-\alpha_{\e}^{p-1-\varepsilon}\lambda_{\e}^{ -\varepsilon (n-2) / 2 }\right) = R_1 ,\notag \\
& 	( EL')  \quad c_2    \e - c(n) \frac{ \ln ^{\s_n} ( \l_{ \e} ) }{ \l_{ \e} ^2 } V(a_{ \e} ) + \ov{c} \frac{H(a,a)}{ \l_\e ^{n-2}} = R_2, \label{EL1} \\   
&(EA')  \quad   c_2(n)  \frac{ \ln ^{\s_n} ( \l_{ \e} ) }{ \l_{ \e} ^3 } \n V(a _{\e} ) - \frac{\ov{c}'}{ \l _\e^{n-1} } \frac{ \partial H}{ \partial a_\e }(a_\e,a_\e)  = R_3, \label{EA1}  
  \end{align}
where $R_1$, $R_2$ and $R_3$ are lower-order terms, and $H$ denotes the regular part of the Green's function of the Laplace operator under Dirichlet boundary conditions.\\
Notice that, for $ n = 3 $, the equivalent of equation \eqref{EL1} is given by 
\be c_2    \e + O \Big( \frac{ 1 }{ \l_{ \e}  } \Big) + \ov{c} \frac{H(a,a)}{ \l _\e} = R_2. \label{EL1n3} \ee

 Observe that
 \begin{itemize}
        \item for $ n = 3$, if the concentration point $ a_\e $ is close to the boundary, then $ H(a_\e, a_\e) $ is very large and it will dominate the term $ O( 1 / \l_\e ) $ and therefore, we are able to get a contradiction. Hence we deduce that, for the dimension $ n=3 $, the concentration point is far away from the boundary, as shown in \cite{FKK}. 
       \item for $ n \geq 4$,  when $V < 0 $, the leading terms in \eqref{EL1} are positive, which leads to a contradiction. This illustrates  that if $V < 0 $, the concentration point remains away from the boundary. Combining this observation with Remark \ref{rn}, we conclude that, for $ V < 0 $, blow-up phenomena do not occur for  problem $(\mathcal{P}_{\e, V})$.
       \item for $ n \geq 4$, if $ V > 0 $, we believe that the concentration point $a_\e$ may converge to a point $ \xi^* \in \partial \O$. In this case, 
 $ \xi^*$ must be a critical point of $V$ restricted to the boundary, i.e., $ V_{ | \partial \O} $, and it must satisfy
$$ \frac{ \partial V}{ \partial \nu }  > 0, $$
where $\nu$ denotes the outward normal to $\partial \O$ at $\xi^*$.
Indeed,   multiplying \eqref{EA1} by $ \nu $ and using the fact that (see \cite{Rg}) 
 $$ \frac{\partial H }{ \partial a_\e }  ( a_\e, a_\e) \cdot \nu = \frac{ c + o(1) }{ d_\e ^{n-1}} \qquad (\mbox{where } \, \, d_\e := d(a_\e , \partial \O) ) , $$
we derive that 
\be \label{nbv2}   \frac{ \ln ^{\s_n} ( \l_{ \e} ) }{ \l_{ \e} ^3 } \frac{ \partial V}{ \partial \xi^*} ( \xi^* ) = \frac{ c + o(1) }{ ( \l _\e d_\e)^{n-1} }  . \ee
This implies that  
$$ \frac{ \partial V}{ \partial \xi^*} ( \xi^* )  > 0  \qquad \mbox{ and } \qquad \frac{ 1 }{ ( \l _\e d_\e)^{n-2} } = ( c + o(1) )  \Big( \frac{ \ln ^{\s_n} ( \l_{ \e} ) }{ \l_{ \e} ^3 } \Big) ^{(n-2) /(n-1) } = o \Big(  \frac{ \ln ^{\s_n} ( \l_{ \e} ) }{ \l_{ \e} ^2 } \Big). $$
Hence, putting this information in \eqref{EL1}, we obtain 
\be\label{nbv1} \frac{ \ln ^{\s_n} ( \l_{ \e} ) }{ \l_{ \e} ^2 }  = \e ( c + o(1)) . \ee
In future work, we will construct solutions $u_\e$ exhibiting clustered multi-bubble behavior, where all concentration points converge to the same point on the boundary. For the single boundary bubble case, the underlying ideas are relatively well understood, and the main challenge lies in obtaining precise expansions. However, in the multi-bubble case, it is essential to carefully analyze the contributions from the interaction terms $\e_{ij}$. Note that the construction of interior multi-bubble solutions has already been carried out: see \cite{2ABE} for isolated bubbles and \cite{Clus} for clustered configurations.
\end{itemize}

 \subsection{ Case of non zero weak limit }

In this subsection, we will assume that $ \o \neq 0 $. First, we will estimate the norm of the function $ v_\e $ in this case. 
\begin{lem} The function $ v _\e $ introduced in \eqref{eq:6} with $ \o \neq 0 $ satisfies
$$ \| v_\e \| \leq    c \,  \begin{cases}
\e +  \sum_{  i } \l_i ^{ (2-n)/2} + \sum_{ j \neq i } \e_{ij} & \mbox{ if } n \leq 5 , \\
\e +  \sum_{  i } \l_i ^{ - 2} \ln ^{2/3} ( \l_i) + \sum_{ j \neq i }  \e_{ij} \ln( \e_{ij}^{-1}) ^{2/3} & \mbox{ if } n = 6 , \\
\e + \sum_{  i } \l_i ^{ - 2} +  \e_{ij}^{(n+2) /( 2(n-2)) } \ln( \e_{ij}^{-1}) ^{(n+2)/(2n)} & \mbox{ if } n \geq 7.
\end{cases} $$
\end{lem}

\begin{pf}
%We will follow the proof of Proposition 2.2 of \cite{Clus}.  
 Multiplying the first equation of the problem $(P_{\e, V})$ by $ v \in E_{\o, a, \l}^\perp $ and integrating over $ \O$, we obtain
\begin{align*} \langle v_\e , v \rangle   = \int_\O u_\e ^{p-\e} v  = &  \int_\O \underline{u}_\e ^{p-\e} v  + \a_0^{p-\e} \int_\O \o ^{p-\e} v + O \Big( \int_\O \underline{u}_\e ^{p-1} | v  | + \int_\O \underline{u}_\e  | v  | \Big) \\
 & + p \sum \int_\O \d_{a_i, \l_i} ^{ p-1} v_\e v + p \int_\O \o^{p-1} v_\e v  + o( \| v_\e \| \| v \| ) 
\end{align*}
which implies that 
\begin{align*} BL_\e (v _\e, v ) & :=  \langle v_\e , v \rangle    - p \sum \int_\O \d_{a_i, \l_i} ^{ p-1} v_\e v  - p \int_\O \o^{p-1} v_\e v + o( \| v_\e \|  \| v \|  ) \\ 
 & = \int_\O \underline{u}_\e ^{p-\e} v  + \a_0 \int_\O \o ^{p-\e} v + \sum O \Big( \int_\O \d_{a_i, \l_i} ^{p-1} | v  | + \int_\O \d_{a_i, \l_i}  | v  | \Big) 
\end{align*}
 Observe that, using Lemma 2.1 of \cite{Clus}, we know that 
$$ \Big|  \int_\O \underline{u}_\e ^{p-\e} v _\e \Big| \leq c \| v _\e \| \Big( \e +   \sum \frac{ \ln^{ \s_n} (\l_i) }{ \l_i ^2 } \Big) + c \| v_\e \|  \sum_{ j \neq i } \begin{cases}
\e_{ij} & \mbox{ if } n \leq 5 , \\
\e_{ij}^{(n+2) /( 2(n-2)) } \ln( \e_{ij}^{-1}) ^{(n+2)/(2n)} & \mbox{ if } n \geq 6.
\end{cases} $$
In addition, since $ v_\e \in E_{\o, a, \l}$, it holds 
$$ \int_\O \o ^{p-\e} v_\e = \int_\O \o ^{p} v_\e + O \Big( \e \int_\O \o ^{p} | v_\e | \Big) = O \Big( \e \| v_\e \| \Big) . $$
Furthermore, it holds
\begin{align*}
 \int_\O \d_{a_i, \l_i} ^{p-1} | v _\e | + \int_\O \d_{a_i, \l_i}  | v _\e | & \leq c \| v _\e \| \Big[ \Big( \int_\O \d_{a_i, \l_i} ^{(p-1) \frac{2n }{n+2}} \Big)^{\frac{n+2}{2n}} + \Big( \int_\O \d_{a_i, \l_i} ^{ \frac{2n }{n+2}} \Big)^{\frac{n+2}{2n}} \Big]  \\
  & \leq c \| v _\e \| \begin{cases} 
1 / \l _i ^{ (n-2) / 2 } & \mbox{ if } n \leq 5 , \\
\ln ^{2/3} ( \l_i) / \l _i ^{ 2 } & \mbox{ if } n = 6 , \\
1 / \l _i ^{ 2 } & \mbox{ if } n \geq 7 . \end{cases}
\end{align*}

Finally, from Lemma $2.3$ of \cite{ABE}, we know that the quadratic form 
$$ Q ( v )  :=  \| v \|^2    - p \sum \int_\O \d_{a_i, \l_i} ^{ p-1}  v ^2  - p \int_\O \o^{p-1}  v^2 $$
is non-degenerate on the space $ E_{\o, a , \l} ^\perp $ which implies that the bilinear form $ BL_\e $ is non-degenerate.  
 Hence the proof of the lemma follows. 
\end{pf}

Regarding the balancing conditions, they become

\begin{pro} \label{balancing om}
For $ \e $ small, the following facts hold 
\begin{align}
& ( EAl_i') \,  \left( 1-\alpha_{i}^{p-1-\varepsilon}\lambda_{i}^{ -\varepsilon \frac{n-2}{ 2 }}\right) = O(R_{\alpha_{i}})  + O\Big(  \frac{1}{ \l_i ^{\frac{n-2}{2}}} \Big)   ,   \label{EAli2}  \\
& 	( EL_i')  \quad c_2    \e - \ov{ c }_2 \sum_{ j \neq i }  \l_{i,\e} \frac{ \partial \e_{ij} }{ \partial \l_{i, \e} } - c(n) \frac{ \ln ^{\s_n} ( \l_{i, \e} ) }{ \l_{i, \e} ^2 } V(a_{i, \e} )\notag \\
&  \qquad \qquad  \qquad = o\Big( \e + \sum_{ j \neq i } \e_{ij}  +\frac{ \ln ^{\s_n} ( \l_{i, \e} ) }{ \l_{i, \e} ^2 } \Big) + O \Big( \| v_\e\|^2 +  \frac{1}{  \l_i ^{\frac{n-2}{2}}}  + \sum_{k \neq j } \e_{kj}^{\frac{n}{n-2}} \ln ( \e_{kj}^{-1}) \Big) , \label{ELi2} \\   
&(EA_i ')  \quad   c_2(n) \a_{i, \e} \frac{ \ln ^{\s_n} ( \l_{i, \e} ) }{ \l_{i, \e} ^3 } \n V(a _{i, \e} ) -  \ov{ c }_2 \sum_{ j \neq i } \a_{j,\e} \frac{1}{\l_{i,\e} }    \frac{ \partial \e_{ij} }{ \partial a_{i, \e} } = O ( R' _{a_i} ) + O\Big( \frac{1}{ \l_i ^{n/2}} \Big), \label{EAi2'}  
  \end{align}
where $R_{\alpha_{i}}$,  $ R'_{a_i} $ and the constants are defined in Proposition \ref{balancing}. 
\end{pro}

\begin{pf} The proof follows the arguments of Propositions $3.3$,  $ 3.4$  and $ 3.5 $ in \cite{Clus}, taking into account the presence of the term $ \a_{0,\e} \o $, and making use of \eqref{om9} along with the following estimate 
\begin{align*} \langle \o , \frac{1}{ \l _i} \frac{ \partial \pi \d_{a_i, \l_i} }{ \partial a_i } \rangle & = \frac{n+2}{n-2} \int_\O \d_{a_i, \l_i}^{4/(n-2)} \frac{1}{ \l_i} \frac{\partial \d_{a_i, \l_i}}{\partial a_i } \o \\
& =  \frac{n+2}{n-2} \o(a_i) \int_{ B_r} \d_{a_i, \l_i}^{\frac{4}{n-2}} \frac{1}{ \l_i} \frac{\partial \d_{a_i, \l_i}}{\partial a_i } + O \Big( \int_{ B_r}  \d_{a_i, \l_i}^{ \frac{n+2}{n-2}} | x - a_i | + \int_{ \R^n \setminus B_r}  \d_{a_i, \l_i}^{ \frac{n+2}{n-2}} \Big) \\
& =  O \Big( \frac{1}{ \l_i ^{ n / 2 }} \Big), 
\end{align*}
where $ B_r := B(a_i, r ) \subset \O$.
 \end{pf}

\begin{rem} \label{rem om2}
(1) Observe that, for $ n \geq 7 $, the remainder term satisfies
$$ \frac{ 1}{ \l_i ^{(n-2) / 2 }} = o \Big( \frac{1}{ \l_i ^2} \Big) $$ 
which implies that the new equations   \eqref{EAli2} and \eqref{ELi2} coincide with equations  \eqref{EAli1} and \eqref{ELi1}. \\
Furthermore, for $ n \geq 8 $, it follows that 
$$   \frac{ 1}{ \l_i ^{n / 2 } } = O \Big( \frac{1}{ \l_i ^4} \Big)  $$
which implies that equation \eqref{EAi2'} coincides with  equation \eqref{EAi1}. 
Thus, the proof done for the zero weak limit will be extended without any change for the dimensions $ n \geq 8$.

However, for $ n = 7 $, the new remainder term $ 1/ \l_i^{n/2} =  1/ \l_i^{7/2}$ does not change the proof of the case  of zero weak limit. But, we notice that the estimate of $ \beta(\e)$ (defined in the theorems) follows from the equation \eqref{tes3} for the zero weak limit. For $ n = 7 $, in this equation, the term $  \frac{ 1}{ \l_i ^{7 / 2 } } $ will be added and it will change the estimate of $ \beta(\e)$. However, for the assertions $(iii)$ of $(b)$ and  $(c)$, the proof remains exactly the same with some minor changes in the power of the $ \e $ in the remainder terms.  

For this reason, in the sequel we will focus on proving Theorem \ref{th:t1}.

(2) For $ n \leq 5 $, we have $ ( n-2) / 2 < 2 $ and therefore the term
$$  1/ \l_i^{( n-2 ) /2} $$
will be the principle term in the equation \eqref{ELi2}. We think that there is no interior blowing up solutions $ u_\e $ which converges weakly to a solution $ \o\neq 0$. This result needs a precise expansions of the equations \eqref{EAli2}, \eqref{ELi2} and \eqref{EAi2'} and a refine analysis of these equations to get a contradiction.  We mention that this phenomena is proved in \cite{Druet2} for a similar problem with nonlinear critical exponent.
\end{rem}

\section{ Estimate of the concentration rates }

In this section, we will prove that the concentration rates are of the same order. To this aim, we begin by outlining the relationship between the interaction of the bubbles and the concentration rates.
\begin{lem} \label{asy4}
$$   \e + \sum_{ k \neq j } \e_{kj} \leq c \sum_{ 1 \leq i \leq N } \frac{ \ln ^{\s_n} ( \l_{i, \e} ) }{ \l_{i, \e} ^2 } \leq c' \, \Big(  \e + \sum_{ k \neq j } \e_{kj}  \Big)   .$$
% Furthermore, ordering the $ \l_i$'s as: $ \l_1 \leq \cdots \leq \l_N$, it holds  $$ \frac{ \l_2}{ \l_1 } \leq c .$$
\end{lem}
\begin{pf} The second inequality follows immediately from $(ELi)$ by using \eqref{estveps}. Concerning the first one, summing $ 2^k (EL_k) $ for $ k \in \{1, \cdots, N\}$ and using \eqref{estveps} and  the fact that 
$$ -  \l_{i,\e} \frac{ \partial \e_{ij} }{ \partial \l_{i, \e} } - 2  \l_{j,\e} \frac{ \partial \e_{ij} }{ \partial \l_{j, \e} } \geq c \e_{ij}  \qquad \mbox{ if  }  \l_i \leq \l_j , $$
the result follows.
 \end{pf}

Ordering the $ \l_i$'s as: $ \l_1 \leq \cdots \leq \l_N$ and using Lemma \ref{asy4}, we deduce that 
\be \label{vep}  \| v_\e \|^2 \leq c  \Big( \frac{ \ln ^{\s_n} ( \l_{1, \e} ) }{ \l_{1, \e} ^2 } \Big) ^2 + c \Big( \sum \e_{k\ell }^{n/(n-2) } \ln ( \e_{k \ell }^{-1} )  \Big) ^{(n+2)/ n } , \qquad  R_{\a_i} \leq c \frac{ \ln ^{\s_n} ( \l_{1, \e} ) }{  \l_{1, \e} ^2 } . \ee
Next, we present the following result, which constitutes the main objective of this section.
\begin{pro} \label{sameorder}  Ordering the $ \l_i$'s as: $ \l_1 \leq \cdots \leq \l_N$, then  there exists a positive constant $ M $ such that 
$$ \l_N  \leq M \l_1 . $$ 
\end{pro}
The sequel of this section is devoted to the proof of Proposition \ref{sameorder}. To make the proof of this proposition clearer, we will split it into several subsidiary results. Arguing by contradiction, we will assume that $ \l_N / \l_1 \to \infty $ and our goal is to get a contradiction. 
Let $$ \mathcal{S}_1:= \{ i \geq 2 : \lim | a_i - a_1 | = 0 \} ,  \quad  \mathcal{S}_2 := \{ i \geq 2 : \lim \l_i / \l_1 = \infty \}  , \quad  i_0 := \min \mathcal{S}_2 .$$
Notice that, by our assumption, we have $ N \in \mathcal{S}_2 $ and therefore $ \mathcal{S}_2 \neq \emptyset $.  The next lemma aims to further clarify the relationship between the interaction of the bubbles and the concentration rates.
\begin{lem} \label{asy50} For each $ j \geq 2 $, it holds 
\be \label{asy7}  c \e + c \sum_{k \geq j} \sum_{ \ell \neq k } \e_{k \ell } \leq c \frac{ \ln ^{\s_n} ( \l_{j } ) }{ \l_{j} ^2 } + c \Big( \frac{ \ln ^{\s_n} ( \l_{1 } ) }{ \l_{1} ^2 } \Big)^{n/(n-2)} \ln ( \l_1)  + \frac{ c }{ ( \l_1 \l_j )^{(n-2)/2}} . \ee 
In addition, taking $ j = i_0$ in the above formula, we derive that 
\be \label{asy51} \e = o \Big( \frac{ \ln ^{\s_n} ( \l_{ 1 } ) }{ \l_{ 1 } ^2 } \Big) . \ee
\end{lem}
\begin{pf} Following the proof of Lemma \ref{asy4}, the first part follows by taking  $ \sum_{k \geq j} 2^k (EL_k) $. The second one follows from the first one and the fact that 
\be\label{asy0}  \frac{ 1  }{ ( \l_1 \l_j )^{(n-2)/2}}  \leq \frac{ 1  }{  \l_1 ^{n-2} } = o  \Big( \frac{ \ln ^{\s_n} ( \l_{ 1 } ) }{ \l_{ 1 } ^2 } \Big) . \ee
Thus the proof is completed. 
\end{pf}

The next result addresses the case where the concentration point is isolated from the other concentration points.
\begin{lem}\label{l34} If $ \mathcal{S}_1= \emptyset $ or $ i_0 = 2 $, then the conclusion of Proposition \ref{sameorder} holds. 
\end{lem}
\begin{pf}
If $ \mathcal{S}_1= \emptyset $, we deduce that $ | a_1 - a_i | \geq c > 0 $ for each $i \geq 2$. Therefore, 
\be \label{asy52} \e_{1 j } \leq  \frac{ c  }{ ( \l_1 \l_j )^{(n-2)/2}}  \leq \frac{ c }{  \l_1 ^{n-2} } = o  \Big( \frac{ \ln ^{\s_n} ( \l_{ 1 } ) }{ \l_{ 1 } ^2 } \Big) \quad \forall \, \, j \geq 2. \ee
Now, if $ i_0 = 2 $, then, taking $ j = i_0 =2$ in Lemma \ref{asy50}, we derive that 
$$ \e_{k \ell }    = o  \Big( \frac{ \ln ^{\s_n} ( \l_{ 1 } ) }{ \l_{ 1 } ^2 } \Big) \quad \forall \, \,  k \neq  \ell . $$
Putting this information in $(EL_1)$ and using \eqref{asy51}, we get a contradiction which ends the proof of Proposition \ref{sameorder}.
\end{pf}

With the help of Lemma \ref{l34}, we can assume that $\mathcal{S}_1 \neq \emptyset $  and  $3 \leq  i_0 \leq N $ in the remainder of this section. Now, let 
\begin{align*} 
 & \vn'_1 := \{ j : \exists \, \, i \in \mS_1, i < i_0  \mbox{ with }  \l_i | a_i -a_j | \nto  \infty \mbox{ as } \e \to 0 \} \quad \mbox{ and } \quad \vn_1 := ( \mS_2 \cap \mS_1)  \setminus \vn' _1, \\
 &  \vn'_2 := \{ j \in \mS_1: \exists \, 1 \leq i \leq N \mbox{ with }  (\l_j^2 / \ln^{\s_n} (\l_j) ) \e_{ij} \to  \infty \mbox{ as } \e \to 0 \} \,  \mbox{ and } \,  \vn_2 := \mS_1 \setminus \vn' _2.
 \end{align*}
It is easy to see  that  $ \vn'_1 \subset \mS_2$. We will now estimate the concentration rates when the indices are in these sets. 
\begin{lem} \label{asy8} $(1)$ Assume that $ \vn'_1 \neq \emptyset$ and  let $ j \in \vn'_1$.  Then, it holds
$$ \frac{1}{\l_j ^{ (n-2)/2} } \leq \frac{ \ln ^{ \s_n} ( \l_1) }{ \l_1 ^{(n+2)/2}} \quad  \mbox{ and } \quad \e + \sum _{ k \geq j; \ell \neq k } \e_{k\ell } \leq c \Big( \frac{ \ln ^{\s_n} ( \l_{1 } ) }{ \l_{1} ^2 } \Big)^{n/(n-2)}  \ln(\l_1) .  $$
$(2)$ Assume that $ \vn'_2 \neq \emptyset$ and  let $ j \in \vn'_2$.  Then, it holds
\be \label{asy54}  \e + \sum _{ k \geq j; \ell \neq k } \e_{k\ell } \leq c \Big( \frac{ \ln ^{\s_n} ( \l_{1 } ) }{ \l_{1} ^2 } \Big)^{n/(n-2)} \ln(\l_1)  . \ee
$(3)$ Let $ j \in \vn_2 \cap \mS_2$.  Then, it holds 
$$ \frac{1}{ \l_i } \Big| \frac{\partial \e_{ij} }{ \partial a_i } \Big| = o \Big( \Big( \frac{  \ln ^{\s_n}  ( \l_1) } { \l_1 ^ 2 } \Big)^{ \frac{n-1} {n-2} } \Big) \quad \forall \, \, i <i_0.$$
\end{lem}
\begin{pf} 
To prove Claim (1), let $ j \in \vn'_1$, from the definition, there exists  $ i \in \mS_1$ with $ i < i_0$ such that $ \l_i | a_i - a_j | $ is bounded.  Thus, it follows that 
$$ \l_j / \l_i \to \infty \qquad  \mbox{ and } \qquad \e_{ij} \geq c ( \l_i / \l_j )^{(n-2)/2} .$$
Therefore,  the first assertion of Claim (1) follows from Lemma \ref{asy4}. Furthermore, the second assertion of Claim (1) follows from the first one and \eqref{asy7}. \\
Concerning Claim (2), it follows from \eqref{asy7} since, for $ j \in \vn'_2 $, the term $ \ln^{\s_n} (\l_j) / \l_j^2 $ (which appears in the right hand side of \eqref{asy7}) is small with respect to at least one of the $ \e_{ij} $'s (which appear in the left hand side). In addition, since $ j \in \mathcal{S}_1$, we derive that $ |a_1-a_j| $ is small and therefore 
$$ \frac{ 1 } { ( \l_1 \l_j ) }  \e_{1j} ^{-2/(n-2) } = \frac{1}{ \l_1^2} +  \frac{1}{ \l_j^2} + | a_1 - a_j |^2 $$
which implies that 
$$ \frac{ 1 } { ( \l_1 \l_j )^{(n-2) / 2 } } = o( \e_{1j} ) . $$
 Thus Claim (2) follows. \\
Finally, concerning Claim (3), let $ j \in \vn_2 \cap \mS_2$  and  $ i < i_0$, observe that, 
\begin{itemize}
\item   if $ \l_j \geq \l_i ^2 $, we see that 
\begin{align*}
 &  \frac{ 1 } { ( \l_1 \l_j )^{(n-2) / 2 } } \leq  \frac{ 1 } {  \l_1 ^{ 3 (n-2) / 2 } } = o \Big( \Big( \frac{  \ln ^{\s_n}  ( \l_1) } { \l_1 ^ 2 } \Big)^{ \frac{n-1} {n-2} } \Big)  ,  \\
 & \frac{  \ln ^{\s_n}  ( \l_j) } {  \l_j ^ 2  } \leq  \frac{  \ln ^{\s_n}  ( \l_1) } {  \l_1 ^{ 4 } } = o \Big( \Big( \frac{  \ln ^{\s_n}  ( \l_1) } { \l_1 ^ 2 } \Big)^{ \frac{n-1} {n-2} } \Big)  .
 \end{align*}
Thus, by easy computations and using \eqref{asy7}, we get
$$   \frac{ 1}{ \l_i } \Big| \frac{\partial \e_{ij} }{ \partial a_i } \Big| \leq  c \l_j | a_i - a_j | \e_{ij}^{ \frac{n }{n-2}}  \leq c \e_{ij}   = o \Big( \Big( \frac{  \ln ^{\s_n}  ( \l_1) } { \l_1 ^ 2 } \Big)^{ \frac{n-1} {n-2} } \Big)  . $$
\item If  $ \l_j \leq \l_i ^2 $, recall that $ i < i_0$ and $ j \in   \vn_2$ with  $ j \geq i_0$, thus it follows that  
\begin{align*} 
\frac{ 1}{ \l_i } \Big| \frac{\partial \e_{ij} }{ \partial a_i } \Big|   \leq c \sqrt{ \frac{ \l_j }{ \l_i} }  \e_{ij}^{ \frac{n-1}{n-2}}  & = c \Big( \frac{\l_j^2 }{ \ln^{\s_n} (\l_j) } \e_{ij} \Big) ^{ \frac{n-1}{n-2}} \Big( \frac{ \l_i }{ \l_j } \Big) ^{ 2 \frac{n-1}{n-2} - \frac{1}{2} } \Big( \frac{ 1 }{ \l_i } \Big) ^{ 2 \frac{n-1}{n-2} } \Big(\ln^{\s_n} (\l_j) \Big) ^{ \frac{n-1}{n-2}}  \\
 & = o \Big( \Big( \frac{  \ln ^{\s_n}  ( \l_1) } { \l_1 ^ 2 } \Big)^{ \frac{n-1} {n-2} } \Big).  
\end{align*}
\end{itemize}
 This completes the proof of Lemma \ref{asy8}. 
\end{pf}

We note that, to complete the proof of Proposition \ref{sameorder}, it remains to consider the case where $ \mS_1 \neq \emptyset$ and $ i_0 \geq 3 $. To this aim, we need to prove the following:
\begin{lem}\label{barycenter}
Let $ \l_1, \cdots , \l_\ell $ be such that $ \l_1 \leq \cdots \leq \l_\ell \leq M \l_1 $ for some positive constant $M$ and $a_1, \cdots , a_\ell \in \O $ be such that $ \e_{ij }$ is small for each $ i \neq j$. Let $ b \in \O$ and $ \a_i >0$ for $ i \in \{ 1 , \cdots , \ell \}$,  it holds
$$ - \sum_{ 1 \leq i \neq j \leq \ell } \a_i \a_j  \frac{ \partial \e_{ij} } { \partial  a_i} (a_i - b ) = (n-2) \sum_{ 1 \leq i < j \leq \ell } \a_i \a_j  \e_{ij} (1+o(1)) . $$
\end{lem}
\begin{pf}
Observe that 
\be \label{eij} \frac{ \partial \e_{ij} } { \partial  a_i} = (n-2) \e_{ij} ^{n/(n-2)} \l_i \l_j (a_j -a_i) = -  \frac{ \partial \e_{ij} } { \partial  a_j} .  \ee
 Furthermore, since $ \l_i $ and $ \l_j$ are of the same order, we get 
 $$ \e_{ij} = \frac{1+ o(1) }{ ( \l_i \l_j | a_i - a_j |^2 )^{ (n-2) / 2 } } . $$
 These imply that 
$$ - \frac{ \partial \e_{ij} } { \partial  a_i} (a_i -b) -  \frac{ \partial \e_{ij} } { \partial  a_j} (a_j -b) =  \frac{ \partial \e_{ij} } { \partial  a_i} (a_j -a_i ) = (n-2) \e_{ij} ^\frac{n}{n-2} \l_i \l_j | a_j -a_i |^2 = (n-2)  \e_{ij}  (1+o(1)) . $$
This completes the proof of the lemma.
\end{pf}

%%%%%%%%%%%%%%%%%%%%%%%%%%%%%%%%%%%%%%%%%%%%%%%%%%%%%%
\begin{pfn}{\bf of Proposition \ref{sameorder} in case where $ \mS_1 \neq \emptyset$ and $ i_0 \geq 3 $.}
%%%%
Observe tha, by easy computations, we have 
$$ - \l_i \frac{\partial \e_{ij}}{\partial \l_i }  = \frac{n-2}{2} \e_{ij} (1+ o(1)) \qquad \forall \, \, i \neq j < i_0 . $$
Thus, summing $(EL_k)$ for $ k \in \mS_3:= \mS_1 \cap \{ 1, \cdots, i_0-1\}$ and using \eqref{asy51} and Lemma \ref{asy4}, we obtain
$$  \ov{c}_2 (n-2)  \sum_{k \neq \ell \in \mS_3 } \e_{k \ell }  - c(n) \sum_{k \in \mS_3} \frac{ \ln ^{\s_n}\l_k}{ \l_k ^2} V(a_k) = o \Big(  \frac{  \ln ^{\s_n}  ( \l_1) } { \l_1 ^ 2 }   \Big) + O \Big( \sum_{ k \in \mS_3 ; \, j \notin \mS_3} \e_{kj} \Big)  . $$
Now, we claim that 
$$ \e_{k \ell } = o \Big(  \frac{  \ln ^{\s_n}  ( \l_1) } { \l_1 ^ 2 }   \Big) \qquad \forall \, \, k \in \mS_3, \, \, \forall \, \, \ell \notin \mS_3 . $$
Indeed, for $ k \in \mS_3$  and  $ \ell \notin \mS_3 $
\begin{itemize}
\item if $ \ell \notin \mS_1$, the claim follows since $ | a_\ell -a_k | \geq c > 0 $ in this case,
\item if $ \ell \in \mS_1 \cap \mS_2 $, then the claim follows from \eqref{asy7} (by taking $j=i_0$) and \eqref{asy0}.
\end{itemize}
Hence the previous formula becomes
\be \label{asy111}
 \ov{c}_2 (n-2)  \sum_{k \neq \ell \in \mS_3 } \e_{k \ell }  - c(n) \sum_{k \in \mS_3} \frac{ \ln ^{\s_n}\l_k}{ \l_k ^2} V(a_k) = o \Big(  \frac{  \ln ^{\s_n}  ( \l_1) } { \l_1 ^ 2 }   \Big) . \ee 

Let $i_1, i_2$ be defined by: $$ | a_{i_1} - a_{i_2} | = \min \{ | a_i - a_j | : i \neq j  \in \mS_3 \} . $$
Thus, we derive that
\be \label{asy112}  \e_{ij} \leq \frac{ c }{ ( \l_i \l_j | a_i - a_j |^2 )^{(n-2) / 2 } } \leq \frac{ c }{ ( \l_{i_1} \l_{i_2} | a_{i_1} - a_{i_2} |^2 )^{(n-2) / 2 } } \leq   c  \e_{ i_1 i_2} \qquad \forall \, \, i, j \in \mS_3 . \ee
Therefore,  \eqref{asy111} and \eqref{asy112} imply that 
\be \label{cafe}  c' \frac{ \ln^{\s_n} ( \l_1) } { \l_1 ^2 } \leq \e_{i_1 i_2 } \leq c \frac{ \ln^{\s_n} ( \l_1) } { \l_1 ^2 } \ee
which implies that 
\be \label{asy113}  \frac{c}{ \ln ^{\s_n} (\l_1)} \frac{1}{ \l_1 ^ { n-4 } } \leq  | a_{i_1} - a_{i_2} | ^{n-2} \leq  \frac{ c' }{ \ln ^{\s_n} (\l_1)} \frac{1}{ \l_1 ^ { n-4 } } . \ee

Let $$ \mS'_4 := \{ j \in \mS_3 : \lim | a_j - a_{i_1} | / | a_{i_1} - a_{i_2} |  = \infty \}  \quad ; \quad \mS_4 := \mS_3 \setminus \mS'_4 . $$
Notice that $ i_1, i_2 \in \mS_4$ and 
\begin{align}
 & | a_k -a_\ell | \leq c | a_{i_1} - a_{i_2} | \quad \forall \, \, k, \ell \in \mS _4 , \label{cafe3}\\
  &  \lim | a_k -a_j | / | a_\ell - a_{i_1}  | = \infty  \quad \forall \, \, k, \ell \in \mS _4 \mbox{ and } \forall \, \, j \notin \mS_4 . \label{cafe2}
  \end{align}
Summing $ \a_k \l_k ( a_k - a_{i_1} ) ( EA_k) $ for $ k \in \mS_4 $,  we get
\begin{align} \label{asy115}  O \Big( \frac{ \ln ^{\s_n}( \l_1)}{ \l_1^ 2 }  & \sum_{ k \in \mS_4} | a_k - a_{i_1} | \Big) - \ov{c}_2 \sum_{ k\neq j  \in \mS_4 } \a_k \a_j \frac{ \partial \e_{kj}} { \partial a_k } (a_k - a_{i_1} ) + O \Big( \sum_{ k \in \mS_4; \, j \notin \mS_4} | a_k - a_{i_1} | \Big| \frac{ \partial \e_{kj}} { \partial a_k } \Big| \Big) \notag \\
 & = \sum_{ k \in \mS_4} O \Big( \l_k | a_k - a_{i_1} |  \Big( R_{ a_{k} } + R_{ \a_{k} } \Big[  \frac{ \ln ^{\s_n} ( \l_{k} ) }{ \l_{k} ^3 } + \sum _{j \neq k}  \frac{1}{\l_{k} }  \Big|  \frac{ \partial \e_{jk} }{ \partial a_{k} } \Big| \Big] \Big)  \Big).  \end{align}
Observe that, Lemma \ref{barycenter}  and  \eqref{cafe}  imply that 
\be \label{3wxc} -  \sum_{ k\neq j  \in \mS_4 } \a_k \a_j \frac{ \partial \e_{kj}} { \partial a_k } (a_k - a_{i_1} ) = (n-2)  \sum_{ k\neq j  \in \mS_4 } \a_k \a_j  \e_{kj} (1+ o(1) ) \geq c \frac{ \ln ^{\s_n}( \l_1)}{ \l_1^ 2 }  .  \ee 
Furthermore, for each $ k \in \mS_4$, we have $ \lim | a_k - a_{i_1} | = 0 $ and therefore  
$$ \frac{ \ln ^{\s_n}( \l_1)}{ \l_1^ 2 } \sum_{ k \in \mS_4} | a_k - a_{i_1} | = o \Big( \frac{ \ln ^{\s_n}( \l_1)}{ \l_1^ 2 }  \Big) . $$
Now, we claim that 
\be \label{asy114}
 | a_k - a_{i_1} | \Big| \frac{ \partial \e_{kj}} { \partial a_k } \Big| = o \Big( \frac{ \ln ^{\s_n}( \l_1)}{ \l_1^ 2 }  \Big)  \qquad \forall \, k \in \mS_4, \, \, \forall \, \, j \notin \mS_4 . 
\ee
Indeed, for $ j \notin \mS_4$, four cases may occur:
\begin{itemize}
\item if $ j \notin \mS_1$, it follows that $ | a_k - a_j | \geq c > 0 $ and therefore the claim follows easily (as in \eqref{asy0}), 
\item if $ j \in \vn'_2 $, using \eqref{cafe2} and  Lemma \ref{asy8}, we get 
$$  | a_k - a_{i_1} | \Big| \frac{ \partial \e_{kj}} { \partial a_k } \Big| = o \Big( \l_k \l_j | a_k - a_j |^2 \e_{kj}^{n/(n-2) } \Big) = o ( \e_{kj} ) = o \Big( \frac{ \ln ^{\s_n}( \l_1)}{ \l_1^ 2 }  \Big) $$ 
which implies the claim \eqref{asy114} in this case.
\item if $ j \in \vn_2 \cap \mS_2 $,  using \eqref{cafe3}, \eqref{asy113} and Lemma \ref{asy8}, we obtain 
$$  | a_k - a_{i_1} | \Big| \frac{ \partial \e_{kj}} { \partial a_k } \Big| = o \Big( \l_k  | a_k - a_{i_1} | \Big[ \frac{ \ln ^{\s_n}( \l_1)}{ \l_1^ 2 }  \Big]^{ (n-1) / (n-2 )} \Big) = o \Big( \frac{ \ln ^{\s_n}( \l_1)}{ \l_1^ 2 }  \Big) $$
which implies the claim \eqref{asy114} in this case.
\item if $ j \in \mS_1 \cap \{ 1, \cdots , i_0-1\}$ and $ j $ satisfies $ \lim | a_j - a_{i_1} | /  | a_{i_2} - a_{i_1} | = \infty $, using Lemma \ref{asy4} and \eqref{cafe2}, we obtain
$$  | a_k - a_{i_1} | \Big| \frac{ \partial \e_{kj}} { \partial a_k } \Big| \leq c  \frac{ | a_k - a_{i_1} | }{| a_k - a_j  | } \e_{kj} \leq c \frac{ | a_{i_2} - a_{i_1} | }{| a_k - a_j  | }  \frac{ \ln ^{\s_n}( \l_1)}{ \l_1^ 2 } = o \Big(  \frac{ \ln ^{\s_n}( \l_1)}{ \l_1^ 2 } \Big) $$
which implies the claim \eqref{asy114}  in this case also. 
\end{itemize}
Thus Claim \eqref{asy114} is thereby proved. \\ 
Finally, using \eqref{asy113}, we get 
$$ \l_k | a_k - a_{i_1} | \leq c \Big(\frac{\l_1^2}{\ln ^{\s_n} \l_1}\Big)^{1/(n-2)} .$$
Thus, using Lemma \ref{asy4} and \eqref{vep}, for $ k \in \mathcal{S}_4$, we derive that 
$$ \l_k | a_k - a_{i_1} | R_{\a_k} \Big[  \frac{ \ln ^{\s_n} ( \l_{k} ) }{ \l_{k} ^3 } + \sum _{j \neq k}  \frac{1}{\l_{k} }  \Big|  \frac{ \partial \e_{jk} }{ \partial a_{k} } \Big| \Big] \leq c \Big( \frac { \l_{1} ^2 } { \ln ^{\s_n} ( \l_{1} ) } \Big)^{1/(n-2)}  \Big(\frac  { \ln ^{\s_n} ( \l_{1} ) } { \l_{1} ^2 }  \Big) \Big[  \frac{ \ln ^{\s_n} ( \l_{1} ) }{ \l_{1} ^2 }  \Big] . $$ 
In addition, using again Lemma  \ref{asy4} and \eqref{asy114}, we get 
$$ \l_k | a_k - a_{i_1} | \e_{j\ell }^{ n/ (n-2) } \ln \e_{j \ell }^{-1} \leq c \Big( \frac { \l_{1} ^2 } { \ln ^{\s_n} ( \l_{1} ) } \Big)^{1/(n-2)} \Big( \frac  { \ln ^{\s_n} ( \l_{1} ) } { \l_{1} ^2 }\Big)^{n/(n-2)} \ln \l_1 = o\Big(  \frac  { \ln ^{\s_n} ( \l_{1} ) } { \l_{1} ^2 } \Big) , $$
$$ \l_k | a_k - a_{i_1} | \l_j | a_k - a_j |  \e_{j\ell }^\frac{ n+1}{n-2 }  \leq c  | a_k - a_{i_1} |  \e_{kj}^\frac{1}{n-2}  \Big( | a_k - a_{i_1} | \Big| \frac{ \partial \e_{kj} }{\partial a_k} \Big|\Big)  = o\Big(  \frac  { \ln ^{\s_n} ( \l_{1} ) } { \l_{1} ^2 } \Big) \, \, \,  \mbox{ for } j \notin \mathcal{S}_4, $$
$$ \l_k | a_k - a_{i_1} | \l_j | a_k - a_j |  \e_{j\ell }^\frac{ n+1}{n-2 }  \leq c \e_{kj}^{ (n-1)/(n-2) } \, \, \,  \mbox{ for } j \in \mathcal{S}_4 . $$
Hence, we derive that
\be \label{2wxc} \l_k | a_k - a_{i_1} |  \Big( R_{ a_{k} } + R_{ \a_{k} } \Big[  \frac{ \ln ^{\s_n} ( \l_{k} ) }{ \l_{k} ^3 } + \sum _{j \neq k}  \frac{1}{\l_{k} }  \Big|  \frac{ \partial \e_{jk} }{ \partial a_{k} } \Big| \Big] \Big)  = o \Big(  \frac{ \ln ^{\s_n}( \l_1)}{ \l_1^ 2 } \Big) \quad \forall \, \, k \in \mS_4 .\ee 
Combining the  estimates \eqref{asy115}-\eqref{2wxc}, we derive a contradiction.\\
Thus the proof of Proposition \ref{sameorder} is completed.
\end{pfn}

%%%%%%%%%%%%%%%%%%%
\section{ Proof of Theorem \ref{th:t1} }
%%%%%%%%%%%%%%%%%%

To clarify the proof of our result, we will begin with two particular cases. 

%%%%%%%%%%%%%%%%%%%%%%%%%%
\subsection{ Proof of Theorem \ref{th:t1} in case where $ | a_i - a_j | \geq c > 0 $ for each $ i \neq j $.}
%%%%%%%%%%%%%%%%%%%%%%%%%%%

First, note that, from Proposition \ref{sameorder}, we know that all the $ \l_i$'s are of the same order. As before, we will order the $ \l_i$'s as: $ \l_1 \leq \cdots \leq \l_N$. Second, in this case, it follows that 
\be \label{po1} \e_{ij} \leq \frac{ c }{ (\l_i \l_j)^{(n-2)/2} } \leq \frac{ c }{ \l_1 ^{n-2} } = o \Big(  \frac{ \ln ^{\s_n}( \l_1)}{ \l_1^ 2 } \Big)  \qquad \forall \, \, i \neq j . \ee
Thus, putting this information in $ (EL_i)$ and using the fact that the $ \l_j$'s are of the same order, we derive that
$$ c_2 \e - c(n)  \frac{ \ln ^{\s_n}( \l_i)}{ \l_i^ 2 } V(a_i) = o \Big(  \frac{ \ln ^{\s_n}( \l_1)}{ \l_1^ 2 } \Big)  \quad \forall \, \, 1 \leq i \leq N  $$
which implies that 
\be \label{RES1} 
\frac{ \ln ^{\s_n}( \l_i)}{ \l_i^ 2 } = \frac{ 1 } {\kappa_1(n) } \frac{ 1}{ V(a_i) } \e ( 1+ o(1)) \qquad \mbox{ where } \quad \kappa_1(n) :=  \frac{ c (n) } {c _ 2 }. \ee
Now, we need to precise the location of the points $ a_i $'s. To this aim, observe that
$$ \frac{ 1 }{ \l_i } \frac{ \partial \e_{ij}}{ \partial a_i } = O \Big( \frac{ 1}{ \l_1 ^{ n-1 }} \Big)  . $$
Thus, putting this information in $(EA_i)$, we obtain 
\be \label{tes3} c_2(n) \a_i  \frac{ \ln ^{\s_n}( \l_i)}{ \l_i^ 3 }  \n V(a_i)  = O \Big(  \frac{ 1 }{ \l_i^ 4 } +  \frac{ 1 }{ \l_i^ {n-1} }  +  \frac{ \ln \l _i}{ \l_i^ 4 } \,  (\mbox{if } n=5) \Big) \ee 
which implies that the concentration point $a_i$ has to converge to a critical point $ y_{j_i} $ of $ V $. 
Finally, since we assumed that the critical points of $ V $ are non-degenerate, we deduce that $ | \n V(a_i ) | \geq c | a_i - y_{j_i} | $. 
This achieves  the proof of the theorem in this case. 

%%%%%%%%%%%%%%%%%%%%%%%%%%
\subsection{ Proof of Theorem \ref{th:t1} in case where $ \lim | a_i - a_j | = 0 $ for each $ i \neq j $.}
%%%%%%%%%%%%%%%%%%%%%%%%%%%
In this case, all the points $ a_i$'s converge to the same point, denoted by $y$.\\
As before, we order the $ \l_i $'s as: $ \l_1 \leq \cdots \leq \l_N$. 
Let $i_1$ and $i_2$ be defined as 
$$ | a_{i_1} - a_{i_2} | := \min \{ | a_i - a_j | : i \neq j \} . $$ 
Since all the $ \l_i $'s are of the same order, there exists  a positive constant $ c > 0$, such that  
\be \label{po3} \e_{ij} \leq c \e_{i_1 i_2 } \qquad \forall \, \, i \neq j . \ee
Furthermore, using Lemma \ref{asy4} and \eqref{vep}, it is easy to deduce that, for each $ i$, it holds
\be \label{Rai2}
R_{a_i} \leq c \, \frac{  \e_{i_1 i_2} }{ \l_1} + c \, \e_{i_1 i_2} ^{ n  / (n-2) } \ln(\e_{i_1 i_2} ^{-1}) +  \Big( \frac{ c }{\l_1 ^3 } (\mbox{if } n=4)\, \,  ; \, \, c \, \frac{ \ln \l_1 }{\l_1 ^4 } (\mbox{if } n=5)  \, \, ; \, \,  \frac{ c }{\l_1 ^4 } (\mbox{if } n \geq 6) \Big) . \ee
Let $$ \mS'_1 := \{ i : \lim ( | a_i - a_{i_1} | / | a_{i_1} - a_{i_2} | ) = \infty \} \qquad \mbox{ and } \qquad  \mS_1 := \{1, \cdots,N \} \setminus \mS'_1. $$
Note that $i_1 , i_2 \in \mS_1$ and for each $i, j \in \mS_1$, we have  $ | a_i -a_j | \leq c  | a_{i_1} - a_{i_2} | $.  Furthermore, we have 
\be \label{j1} \Big| \frac{ \partial \e_{ij} }{\partial a_i } \Big| | a_i - a_{i_1} |  \leq c \e_{ij} \frac{  | a_i - a_{i_1} |  }{ | a_i - a_j | } = o( \e_{ij} ) = o( \e_{ i_1 i_2} )  \qquad \forall \, \, i \in \mS_1 \mbox{ and } \forall \, \, j \notin \mS_1 . \ee

We start by the following lemmas:
\begin{lem} \label{po2}
$$ \lim  \ln^{\s_n} (\l_1) \l_1 ^{n-4}    | a_{i_1} - a_{i_2} |^{n-2}  = \infty .$$
\end{lem}
\begin{pf}
Arguing by contradiction, assume that there exists $ \beta > 0 $ such that 
$$ \ln^{\s_n} (\l_1) \l_1 ^{n-4}  \leq \beta \frac{ 1 } { | a_{i_1} - a_{i_2} |^{n-2}}  $$
which implies that 
\be \label{po6}   \frac{ \ln^{\s_n} (\l_1) }{  \l_1 ^ 2 }   \leq \beta \frac{ 1 } { \l_1 ^{n-2}  | a_{i_1} - a_{i_2} |^{n-2}}   \leq c \e_{i_1 i_2} . \ee
Now, summing $\a_i \l_i (a_i - a_{i_1} ) (EA_i)$ for $ i \in \mS_1$ (we will write $ \sum_{ j \neq i} = \sum_{ j \neq i ; \, j \in \mS_1} + \sum_{ j \notin \mS_1} $)  and using \eqref{Rai2}, \eqref{j1} and  Lemma \ref{barycenter},  we get 
$$ \sum_{ i \neq j \in \mS_1} \e_{ij}   \leq c \sum_{ i \in \mS_1} \frac{ \ln^{\s_n} (\l_1) }{  \l_1 ^ 2 }  | \n V(a_i ) | | a_i - a_{i_1} | + c \l_1 \sum_{ i \in \mS_I} | a_i - a_{i_1} | R'_{a_i}  =  o \Big( \frac{ \ln^{\s_n} (\l_1) }{  \l_1 ^ 2 } + \e_{i_1 i_2}  \Big)  $$ 
which is a contradiction with \eqref{po6}. Hence, the proof of the lemma is completed. 
\end{pf}
\begin{lem} \label{po4}  For each $ i \in \{ 1, \cdots, N \}$, it holds 
$$ \frac{ \ln^{\s_n} (\l_i) }{  \l_i ^ 2 } = \frac{ 1 }{ \kappa_1(n)} \frac{1}{ V(y) } \e (1+ o(1)) $$
where $y$ is the limit of the points $a_i$'s and $ \kappa_1(n) $ is defined in \eqref{RES1}.
\end{lem}
\begin{pf}
Using Lemma \ref{po2} and  Proposition \ref{sameorder}, we deduce that 
$$ \e_{ij} = o\Big( \frac{ \ln^{\s_n} (\l_1) }{  \l_1 ^ 2 }  \Big)  \qquad \forall \, \, i\neq j \in \{ 1, \cdots, N \} . $$
Putting this information in equation $ (EL_i)$, the result follows. 
\end{pf}
\begin{lem} \label{po5}
The limit  $ y $ of the points $ a_i$'s is a critical point of $ V $.
\end{lem}
\begin{pf} We distinguish four cases and in each one we will prove the result. 
\begin{itemize}
\item If there exist $ i \neq j $ such that $ \lim | \n V(a_i ) | /  | \n V(a_j) | = 0 . $ Since $ | \n V| $ is bounded, we deduce that $\lim  | \n V(a_i) | = 0 $ and therefore $ | \n V(y) | = 0$ which gives the result in this case. 
\item If there exists $ i $ such that 
\be \label{j2}  \lim \frac{ \ln ^{\s_n}( \l_i)}{ \l_i^ 3 } |  \n V(a_i) | \e_{i_1i_2}^{-(n-1)/(n-2)} = \infty . \ee
Notice that, \eqref{j2} implies that
$$ \frac{1}{ \l_i} \Big| \frac{ \partial \e_{ij} }{ \partial a_i } \Big| \leq \frac{ c }{ \l_i | a_i - a_j | } \e_{ij} \leq  c \e_{i_1i_2} ^{(n-1) / (n-2 ) } = o \Big( \frac{ \ln ^{\s_n}( \l_i)}{ \l_i^ 3 } |  \n V(a_i) | \Big) . $$
Furthermore, easy computations imply that 
\be\label{tes1} R' _{a_i} = o \Big( \frac{ \ln ^{\s_n}( \l_1)}{ \l_1^ 3 } \Big) . \ee
Thus, using $(EA_i)$, we derive that 
 $$ \frac{ \ln ^{\s_n}( \l_i)}{ \l_i^ 3 } |  \n V(a_i) | = o \Big(  \frac{ \ln ^{\s_n}( \l_i)}{ \l_i^ 3 } \Big) $$ which implies the result. 
\item If for each $ i $, $$ \lim \frac{ \ln ^{\s_n}( \l_i)}{ \l_i^ 3 } |  \n V(a_i) | \e_{i_1i_2}^{-(n-1)/(n-2)} = 0 .$$
Summing $\a_i \l_i (a_i - a_{i_1})(EA_i)$ for $ i \in \mS_1$ (we will write $ \sum_{ j \neq i} = \sum_{ j \neq i ; \, j \in \mS_1} + \sum_{ j \notin \mS_1} $) and using Lemma \ref{barycenter} and \eqref{j1}, we get 
$$ \sum_{i\neq j \in \mS_1} \e_{ij} (1+o(1)) = \sum_{i \in \mS_1} O \Big(  \frac{ \ln ^{\s_n}( \l_i)}{ \l_i^ 2 } |  \n V(a_i) | | a_i - a_{i_1} | + \l_i | a_i - a_{i_1} | R'_{a_i} \Big) .$$ 
Observe that, for $ i \in \mS_1$, it follows that $ | a_i - a_{i_1} | $ and $ | a_{i_2} - a_{i_1} | $ are of the same order which implies that $ \l_i | a_i - a_j | $ and $ \e_{i_1i_2}^{-1 / (n-2) } $ are of the same order. Hence, in this case, we obtain
$$  \frac{ \ln ^{\s_n}( \l_i)}{ \l_i^ 2 } |  \n V(a_i) | | a_i - a_{i_1} |  = o \Big( \e_{i_1i_2}^{\frac{n-1}{n-2}}  \l_i | a_i - a_{i_1} | \Big) = o ( \e_{i_1i_2}),   $$
and easy computations imply that 
$$ \l_i | a_i - a_{i_1} |  R'_{a_i} = o (  \e_{i_1i_2} ) + O\Big( \frac{ | a_{i_2} - a_{i_1} | }{ \l_1^2} \, \, (\mbox{if } n=4) ; \quad \frac{ \ln ( \l_1)}{ \l_1^ 3 } | a_{i_2} - a_{i_1} |  \, \, (\mbox{if } n \geq 5) \Big) \quad \forall i \in \mS_1 .$$
Hence we derive that 
$$ \e_{i_1i_2}  \leq c  \Big( \frac{ | a_{i_2} - a_{i_1} | }{ \l_1^2} \, \, (\mbox{if } n=4) ; \quad \frac{ \ln ( \l_1)}{ \l_1^ 3 } | a_{i_2} - a_{i_1} |  \, \, (\mbox{if } n \geq 5) \Big) . $$
Now, using the fact $ \frac{ \ln ^{\s_n}( \l_i)}{ \l_i^ 3 } |  \n V(a_i) | \leq  \e_{i_1i_2}^{(n-1)/(n-2)} $, we derive that $ |  \n V(a_i) | = o(1)$ which implies the result in this case.
\item It remains the case where: there exist two constants $ \beta_1$ and $ \beta _2$ such that  
\be \label{po9} \beta_1  \e_{i_1i_2}^{(n-1)/(n-2)}  \leq  \frac{ \ln ^{\s_n}( \l_i)}{ \l_i^ 3 } |  \n V(a_i) | \leq \beta_2  \e_{i_1i_2}^{(n-1)/(n-2)} \qquad \forall i . \ee
Summing $\a_i \l_i (EA_i)$ for $ i \in \{1, \cdots, N \}$ and using \eqref{eij} and \eqref{po9}, we get 
$$ c_2(n) \sum_{i=1}^N \a_i ^2  \frac{ \ln ^{\s_n}( \l_i)}{ \l_i^ 2 } \n V(a_i) = \sum_{i=1}^N O \Big( \l_i R'_{a_i}  \Big) = o\Big( \frac{ \ln^{\s_n} (\l_1) }{  \l_1 ^ 2 }  \Big)  . $$
 Using the fact that $ \n V(a_i)= \n V(y) + o(1) $, the result follows in this case also.
\end{itemize}
Hence the proof of the lemma is completed.
\end{pf}

To proceed further, notice that Lemma \ref{po4} implies that 
\be \label{j5} \frac{ c + o(1) }{ \l_i } = \Big( \sqrt{ \e / | \ln \e | } \, \, (\mbox{if } n = 4) \, \, ; \quad \sqrt{ \e } \, \, (\mbox{if } n \geq 5) \Big) . \ee 
In addition, in the following, using \eqref{j5},   Lemma  \ref{asy4} and Proposition \ref{sameorder}, we deduce that 
\begin{align}
& \a_i = 1 + \sum O \Big( \frac{ \ln ^{ \s_n + 1 } ( \l_i) } { \l_i ^2 } \Big) = 1 + O ( \e | \ln \e | ) ,  \label{ALPHA2} \\
& R_{a_i} ' =   o\Big( \e_{i_1i_2} ^{ (n-1) / (n-2) } \Big) + \sum O \Big( \frac{1} { \l_i ^3} \, \, (\mbox{if } n = 4) \, ; \quad   \frac{ \ln \l_i } { \l_i ^4} \, \, (\mbox{if } n = 5) \, ; \quad\frac{1} { \l_i ^4} \, \, (\mbox{if } n \geq 6) \Big) \leq R'',  \label{RA2} \\
 & R''=   o\Big( \e_{i_1i_2} ^{ (n-1) / (n-2) } \Big) +  O \Big( ( \e / | \ln \e | )^{ 3/2} \, \, (\mbox{if } n = 4) \, ; \quad   \e ^2 | \ln \e |  \, \, (\mbox{if } n = 5) \, ; \quad \e^2 \, \, (\mbox{if } n \geq 6) \Big) .       \label{RA22}
\end{align}
Hence, the equation $(EA_i) $, defined in \eqref{EAi1}, becomes 
\be \label{EAi2}
(EA_i ')  \quad  c_2(n)  \frac{ \ln ^{\s_n} ( \l_{i, \e} ) }{ \l_{i, \e} ^2 } \n V(a _{i, \e} ) -  \ov{ c }_2 \sum_{ j \neq i }   \frac{ \partial \e_{ij} }{ \partial a_{i, \e} } = O ( \l_i R _{a_i} )  = O ( \l_i R'' ) \quad \forall \, \, i \in \{ 1, \cdots, N \} .
\ee

%%%%%%%%%%%%%%%%%%%%%%%%%%%%
\begin{pro} \label{pro11}  Let $ n \geq 4 $, there exists a positive constant $ \ov{ c } > 0 $ such that 
\begin{align*} 
& (1) \qquad | \ln \e |^{\s_n/4} \e^{ - (n-4) / (2n) } | a_i - a_j | \geq  \ov{c } ^{-1}  \qquad \forall \, \, i \neq j , \\
& (2) \qquad  | \ln \e |^{\s_n/4} \e ^{ - (n-4) / (2n) } | a_i - y | \leq \ov{ c } \qquad \forall \, \, i \in \{ 1, \cdots, N \} , \\
& (3) \qquad | \ln \e |^{\s_n/4} \e ^{ - (n-4) / (2n) } | a_i - a_j | \leq \ov{ c } \qquad \forall \, \, i \neq j ,
\end{align*}
where $ \s_4 = 1 $ and $ \s_n = 0 $ if $ n \geq 5$.
\end{pro}
\begin{pf}
Without loss of generality, we can assume that 
\be \label{a1a2}  | a_1 - a_2 | = \min \{ | a_i - a_j | : \, i \neq j \} .\ee
Combining Proposition \ref{sameorder} and \eqref{a1a2}, we deduce that 
 \be \label{eij12}
 \e_{ij} =  \frac{ 1+ o(1) } { (\l_i \l_j | a_i - a_j |^2 )^{ (n-2) / 2 } }  \leq c \,  \e_{12} \qquad \forall \, \, i \neq j  . 
 \ee
Let $ i _0 $ be such that $ | \n V(a_{i_0} ) | := \max \{ | \n V(a_i) | : 1 \leq i \leq N \} $. We start by proving the following claim:
 \be \label{claim1} \, \,  \mbox{ There exists }  \, \, M > 0 \, \,  \mbox{ such that }  \, \, \frac{ \ln ^{\s_n} ( \l_{i_0}) }{ \l_{i_0} ^3}  | a_{i_0} - y | \leq M \, \e_{12}^{ (n-1)/(n-2)} . \ee
 To prove \eqref{claim1}, we argue by contradiction. Assume that such a $ M $ does not exist. Thus we derive that 
 $$ \lim  \frac{ \ln ^{\s_n} ( \l_{i_0}) }{ \l_{i_0} ^3}  | a_{i_0} - y |  \, \e_{12}^{ - (n-1)/(n-2)} = \infty .$$
Thus, using \eqref{eij12}, the fact that $y$ is a non-degenerate critical point of $ V $ (which implies that $ c | a_{i_0} - y | \leq | \n V(a_{i_0}) | $), we get 
\be \label{eij11}  \frac{ 1}{ \l_i } \Big | \frac{ \partial \e_{ij} }{ \partial a_i } \Big | \leq c \,  \e_{ ij } ^{ (n-1)/(n-2)} = o \Big( \frac{ \ln ^{\s_n} \l_{i_0} } { \l_{i_0} ^3} | a_{i_0} - y  | \Big) = o \Big( \frac{ \ln ^{\s_n} \l_{i_0} } { \l_{i_0} ^3} | \n V(a_{i_0} ) | \Big) \qquad \forall \, \,  i \neq j .\ee
Now, using \eqref{EAi2} and \eqref{eij11}, we obtain 
$$   \frac{ \ln ^{\s_n} ( \l_{i_0} ) }{ \l_{i_0} ^3 }  | \n V(a _{i _0} ) |  \leq   c   R_{a_i}   $$
which implies, by using \eqref{RA2} and \eqref{eij11}, that
\be \label{QA12}  \frac{ \ln ^{\s_n} ( \l_{i_0} ) }{ \l_{i_0} ^3 }  | a _{i _0} - y  |  \leq   c  \sum \Big( \frac{1} { \l_i ^3} \, \, (\mbox{if } n = 4) \, ; \quad   \frac{ \ln \l_i } { \l_i ^4} \, \, (\mbox{if } n = 5) \, ; \quad\frac{1} { \l_i ^4} \, \, (\mbox{if } n \geq 6) \Big) .\ee
Notice that, from the definition of $ i_0$, we have $ | a_i - y | \leq c | a_{i_0} -y | $ for each $i$ and therefore, for each $ i \neq j$, it follows that $ | a_i - a_j | \leq  | a_i - y | +  | y - a_j | \leq c  | a_{i_0} - y | $. Thus, using Proposition \ref{sameorder} and \eqref{QA12}, it follows that 
$$ \e_{ij} = \frac{1+o(1) } { ( \l_i \l_j | a_i - a_j | ^2 )^{(n-2) / 2 } } \geq  \Big( c \frac{ \ln^2 ( \l_{i_0} ) } { \l_i \l_j} \, \, (\mbox{if } n = 4) \, ; \quad   \frac{ c } { \ln ^3( \l_{i_0}) }\, \, (\mbox{if } n = 5) \, ; \quad c  \, \, (\mbox{if } n \geq 6) \Big) 
$$
which contradicts Lemma \ref{asy4}. Hence the proof of \eqref{claim1} follows. \\
Let 
$$ \mathcal{S}'_5 := \{ i : \lim |a_i - a_1 | / | a_1 - a_2 | = \infty \} \qquad \mbox{ and } \qquad  \mathcal{S}_5 = \{ 1 , \cdots , N \} \setminus \mathcal{S}'_5 . $$
Using Proposition \ref{sameorder}, it follows that
\be \label{QA13}
  | a_1 - a_2 | = o( |a_i - a_j | )   \quad \mbox{ and } \quad \e_{ij} = o( \e_{12}) \qquad \forall \, \, i \in \mathcal{S}_5 , \, \, \, \forall \, \, j \in \mathcal{S}'_5.
\ee
Let $ \ov{a}_1 $ be the barycenter of the points $ a_i $'s, for $ i \in \mathcal{S}_5 $, i.e., $ \ov{a} _1 $ satisfies 
\be \label{pbarycenter} \sum_{ i \in \mathcal{S}_5 } ( a_i - \ov{a}_1 ) = 0 .\ee
It is easy to see that 
\be \label{ak1} | a_k - \ov{a}_1 | \leq c \sum_{ i \in \mathcal{S}_4 } | a_i - a_k |   \leq c | a_1 - a_2 | \qquad \forall \, \, k \in   \mathcal{S}_5. \ee
Multiplying \eqref{EAi2} by $ ( a_i - \ov{a}_1 ) $, we obtain 
\be \label{QA14}
 c_2(n) \sum_{ i \in \mathcal{S}_5}  \frac{ \ln ^{\s_n} ( \l_{i} ) }{ \l_{i} ^2 } \n V(a _{i} ) ( a_i - \ov{a}_1 ) -  \ov{ c }_2 \sum_{ i \in \mathcal{S}_5}\sum_{ j \neq i }   \frac{ \partial \e_{ij} }{ \partial a_{i} }( a_i - \ov{a}_1 ) = \sum_{ i \in \mathcal{S}_5} O ( \l_i | a_i - \ov{a}_1  | R'_{a_i} ) 
\ee
where $ R'_{a_i} $ is defined in \eqref{RA2}. 
Note that, using Lemma \ref{barycenter}, we have 
\be \label{j8}
 - \sum_{ i \in \mathcal{S}_5} \sum _{ j \in   \mathcal{S}_5; \, j \neq i}  \frac{ \partial \e_{ij} }{ \partial a_{i} }( a_i - \ov{a}_1 )   = (n-2)  \sum _{ j < i ;  \, i, j \in  \mathcal{S}_5}  \e_{ij} (1+o(1)) . 
\ee
In addition, using \eqref{QA13}, for each $  i \in \mathcal{S}_5 $ and  $ j \in \mathcal{S}'_5 $,  we get 
$$  \Big| \frac{ \partial \e_{ij} }{ \partial a_{i} } \Big| |  a_i - \ov{a}_1  | = (n-2) \l_i \l_j | a_i - a_j | \e_{ij}^{ n / (n-2) } | a_i - \ov{a}_{1} |  = o ( \l_i \l_j | a_i - a_j |^2 \e_{ij}^{ n / (n-2) } ) = o ( \e_{12})  . $$
Now, we need to study the first term in \eqref{QA14}. Recall taht, by the definition of $i_0$, we have $ |a_i - y | \leq c | a_{i_0}-y | $ for each $i$.  
Note that \eqref{claim1} and \eqref{ak1} imply that 
\begin{align}
&  \l_i | a_i - \ov{a}_1  | \e_{12} ^{ (n-1)/(n-2)} \leq c \l_i | a_1 - a_2 |   \e_{12} ^{(n-1) / (n-2) }   \leq c \,  \e_{12}  \qquad \forall \, \, i \in  \mathcal{S}_5 ,  \label{6wxc} \\
& \frac{ \ln ^{ \s_n} \l_i }{ \l_i ^2} | \n V (a_i) | | a_i - \ov{a}_1 | \leq c \l_i | a_i - \ov{a}_1 | \e_{12} ^{(n-1) / (n-2) } \leq c \e_{12} . \label{j7}
\end{align}
Thus, expanding $ V $ around $ \ov{a}_1 $, and using \eqref{barycenter}, \eqref{j7} and Lemma \ref{po4}, it follows that
\begin{align}
 \sum_{ i \in \mathcal{S}_5}  \frac{ \ln ^{\s_n} ( \l_{i} ) }{ \l_{i} ^2 } \n V(a _{i} ) ( a_i - \ov{a}_1 ) & =   \sum_{ i \in \mathcal{S}_5} \frac{ 1}{ \kappa_1(n) V(y) } \e ( 1 + o(1))  \n V(a _{i} ) ( a_i - \ov{a}_1 )  + o ( \e_{12} ) \nonumber  \\
 & = \frac{ \e }{ \kappa_1(n) V(y) }   \n V( \ov{ a}_1) \Big( \sum_{ i \in \mathcal{S}_5} (a_i - \ov{a}_1) \Big) + O ( \e  | a_1 - a_2 |^2) + o ( \e_{12} ) \nonumber  \\
  & = O ( \e  | a_1 - a_2 |^2)     + o ( \e_{12} ) . \label{4wxc}
  \end{align}
It remains to study the right hand side of \eqref{QA14}. Observe that, using \eqref{ak1}, we have
$$ \l_i | a_i - \ov{a}_1  | \e_{12} ^{ (n-1)/(n-2)} \leq c \l_i | a_1 - a_2 |   \e_{12} ^{(n-1) / (n-2) }   \leq c \,  \e_{12}  \qquad \forall \, \, i \in  \mathcal{S}_5 . $$
Furthermore, for $ n = 4 $, we have 
$$  \l_i | a_i - \ov{a}_1  | \frac{1}{ \l_i ^3} = | a_i - \ov{a}_1  | | a_1 - a_2 |^2  \frac{ 1 }{ \l_i ^2 | a_1 - a_2 |^2 } = o ( \e_{12} )  \qquad \forall \, \, i \in  \mathcal{S}_5 . $$
Thus, for $ n = 4 $, equations \eqref{QA14} and \eqref{j8} imply that 
$$ \frac{ c }{ \l_1^2 | a_1 - a_2 |^2} \leq  \e_{12} \leq \sum_{ i\neq j \in \mathcal{S}_5 } \e_{ij}   \leq c  \e  | a_1 - a_2 |^2   \leq c  \frac{\ln(\l_1)}{ \l_1^2}  | a_1 - a_2 |^2 $$ 
which implies that $ | a_1 - a_2 |^4 \geq c / \ln \l_1 $. Thus, Assertion $(1)$ for $ n=4 $ follows. \\
For $ n \geq 6 $, it is easy to get that 
$$  \l_i | a_i - \ov{a}_1  | \frac{1}{ \l_i ^4} \leq  \frac{ c }{ \l_i  | a_1 - {a}_2  |} | a_1 - {a}_2  |^2 \frac{1}{ \l_i ^2} = o ( \e | a_1 - {a}_2  |^2)  \qquad \forall \, \, i \in  \mathcal{S}_5 . $$
Thus, for $ n \geq 6 $, equations \eqref{QA14}, \eqref{j8} and \eqref{4wxc} imply that 
$$ \frac{ c \, \e ^{(n-2) / 2 } }{ | a_1 - a_2 |^{n-2} } \leq \frac{ c  }{ ( \l_1 | a_1 - a_2 | ) ^{n-2} } \leq 
 \e_{12} \leq \sum_{ i\neq j \in \mathcal{S}_5 } \e_{ij}   \leq c  \e  | a_1 - a_2 |^2  $$ 
which implies Assertion $(1)$ for $ n \geq 6 $. \\
Finally, for $ n = 5 $, equations \eqref{QA14} and \eqref{j8} imply that 
\begin{align}
 \frac{ c \, \e ^{ 3 / 2 } }{ | a_1 - a_2 |^{ 3 } } \leq \frac{c } { \l_1 ^{ 3 } | a_1 - a_2 |^{ 3 } } \leq  \e_{12} \leq \sum_{ i\neq j \in \mathcal{S}_5 } \e_{ij}    & \leq c  \frac{1}{ \l_1 ^2 }  | a_1 - a_2 |^2 + c  \frac{ \ln (\l_1 ) }{ \l_1 ^3} | a_1 - a_2 | \notag \\
  & \leq c \, \e \, | a_1 - a_2 |^{ 2 } + c \, \e^{3/2} | \ln \e | | a_1 - a_2 | . \label{nb1}
  \end{align}
  Two cases may occur:
  \begin{itemize}
\item if $  | a_1 - a_2 | \geq \sqrt{ \e } | \ln \e | $. In this case, \eqref{nb1} implies that 
$$  \frac{  \e ^{ 3 / 2 } }{ | a_1 - a_2 |^{ 3 } } \leq c \, \e \, | a_1 - a_2 |^{ 2 } $$ 
which gives the proof of Assertion $(1)$ easily.
\item if $  | a_1 - a_2 | \leq \sqrt{ \e } | \ln \e | $. In this case,  \eqref{nb1} implies that 
$$  \frac{  \e ^{ 3 / 2 } }{ | a_1 - a_2 |^{ 3 } } \leq c \, \e^{3/2}  | \ln \e | | a_1 - a_2 | $$
and therefore 
$$ | a_1 - a_2 | ^4 \geq c | \ln \e | ^{-1} \geq \e ^{ 4/ 10 } $$
which implies the proof of Assertion $(1)$ .
\end{itemize}
Thus, the proof of Assertion $(1)$ is proved in the case $ n = 5 $. \\
Hence Assertion $ (1) $ is completed. \\
Concerning Assertion $ (2) $, it follows from the first one and equation \eqref{claim1} by using the fact that $ | a_i - y | \leq c | a_{i_0} - y |$ for each $i \in \{ 1, \cdots, N \}$ (from the definition of $ i_0$).\\
Finally, Assertion $(3)$ follows from the second one by using the fact fact $ | a_i - a_j | \leq  | a_i - y | + | a_j - y |$ for each $i, j  \in \{ 1, \cdots, N \}$. 
\end{pf}

\begin{pro} \label{charac-points}
Let $ a_i$'s, for $ i \in \{1,\cdots, N \}$, be the concentration points of the solution $u_\e$. For $ i \in \{1,\cdots, N \}$, let us define $z_i$ as 
\begin{align} 
 & a_i := y + \s \eta(\e) z_i  \qquad \mbox{ for }  i \in \{1,\cdots, N \}  \qquad \mbox{ where }  \label{change2} \\ 
 &  \s^{-1} :=   \kappa_2(n) V(y)^{(n-4) / (2n) } \qquad  \mbox{ with  } \qquad \kappa_2(n)  := \Big(  \frac{c_2(n)}{ \ov{c}_2} \Big)^{1/n} \Big( \frac{c(n) } {c_2} \Big)^{(n-4) / (2n) } ,  \label{kappa2} 
\end{align}
$ \eta (\e) $ is defined in \eqref{etaeps} and  the constants $ c_2(n)$, $ c(n) $, $ c_2 $ and $ \ov{c}_2 $ are introduced in Theorem \ref{th:t1}.\\
Then $(z_1, \cdots, z_N)$ converges to  a critical point of the function $ F_{y , N } $ which is defined in \eqref{FYN}.
\end{pro}

\begin{pf} First, we remark that Proposition \ref{pro11} implies that the new variables $ |z_i| $'s are bounded above and the quantities $ | z_i-z_j|$'s, for $ i \neq j$, are bounded below and above. \\
Now, we need to rewrite the equation \eqref{EAi2}. To this aim, expanding $ \n V $ around $y$ and using \eqref{change2}, we get
$$ \n V(a_i) = D^2 V(y) (a_i - y , . ) + O ( | a_i - y |^2) = \s \eta (\e) D^2 V(y) (z_i  , . ) + O (  \eta(\e) ^2 | z_i  |^2) . $$
Furthermore, using again \eqref{change2} and Lemma \ref{po4}, we get
\begin{align*} 
\frac{\partial \e_{ij} }{ \partial a_i } & = (n-2) \l_i \l_j (a_j - a_i ) \e_{ij}^{n/(n-2)} = (n-2) \frac{ \l_i \l_j \s \eta(\e) (z_j-z_i) }{ ( \l_i \l_j \s^2 \eta(\e)^2 | z_j - z_i |^2 )^{n/2} } (1+o(1)) \\
& = \frac{ n-2}{ [\s \eta(\e)] ^{n-1}}  \Big( \frac{ c_2}{c(n)} \frac{ \e }{V(y)}\Big)^{(n-2)/2}  \Big( \frac{ 2 }{ | \ln \e | } \Big)^{\s_n} \frac{ z_j - z_i }{ | z_j - z_i |^ n }  + o\Big( \frac{ \e^{(n-2)/2}}{ \eta(\e)^{n-1} } | \ln \e |^{-\s_n}  \Big) .
\end{align*}
In addition, using \eqref{RA2}, Proposition \ref{pro11} and Lemma \ref{po4},  we see  that 
\be\label{tes2} \l_i R'_{a_i} = o( \e \eta(\e) ) \qquad \forall \, \, i \in \{ 1, \cdots, N \} . \ee  
Thus, putting these informations in \eqref{EAi2} and using Lemma \ref{po4}, we obtain
\begin{align*} c_2(n) \Big( \frac{ c_2}{c(n)} \frac{ \e }{V(y)}\Big) & \s \eta(\e) D^2 V (y) ( z_i , . ) \\
 &   - \ov{c}_2 \frac{ n-2}{ [\s \eta(\e)]^{n-1} }  \Big( \frac{ c_2}{c(n)} \frac{ \e }{V(y)}\Big)^\frac{n-2}{2}  \Big( \frac{ 2 }{ | \ln \e | } \Big)^{\s_n} \sum_{ j \neq i }   \frac{ z_j - z_i }{ | z_j - z_i |^ n } 
  = o ( \e \eta(\e) ) \end{align*}
which implies that (by the choose of the values of $ \s $ and $ \eta (\e) $)
$$  D^2 V (y) ( z_i , . )  - (n-2)  \sum_{ j \neq i }   \frac{ z_j - z_i }{ | z_j - z_i |^ n } = o ( 1 )  \qquad \forall \, \, i \in \{ 1, \cdots, N \}. $$
Taking $ \ov{z}_i := \lim z_i$, we deduce that 
$$  D^2 V (y) ( \ov{z}_i , . )  - (n-2)  \sum_{ j \neq i }   \frac{ \ov{z}_j - \ov{z}_i }{ | \ov{z}_j - \ov{z}_i |^ n } = 0  \qquad \forall \, \, i \in \{ 1, \cdots, N \}. $$
Hence, the point $ (\ov{z}_1, \cdots, \ov{z}_N) $ is a critical point of the function $ F_{y,N} $. This completes the proof of the proposition. 
\end{pf}

\begin{pfn}{\bf of Theorem \ref{th:t1} in case where $ \lim | a_i - a_j | = 0 $ for each $ i \neq j$. }
The proof follows from Propositions \ref{po4}, \ref{pro11} and \ref{charac-points}. 
\end{pfn}

%%%%%%%%%%%%%%%%%%%%%%%%%%%%%%%%%%
\subsection{ Proof of Theorem \ref{th:t1} in the general case }

Notice that we have proved that the concentration rates $ \l_i$'s are of the same order, see Lemma \ref{po4}. Hence, we will denote by $ O(f(\l)) $ each quantity of the form $ \sum O(f(\l_i))$. \\
Let $ i \in \{ 1, \cdots, N \}$, two cases may occurs:
\begin{description}
\item [ (i) ] The first one, if there exists a constant $ c > 0 $ such that $ | a_i - a_j | \geq c $ for each $ j \neq i$. In this case, according to Subsection 4.1, we know that $ a_i$ converges to a critical point $y_{j_i}$. Furthermore, for $ r > 0 $ small, the ball $ B(y_{j_i} , r)$ contains only $a_i$ for $ \e$ small. 
\item[ (ii)] The second case, if there there exists $ j \neq i $ such that $ \lim | a_i - a_j | =0 $. In this case, let us introduce the set 
$$ A^ i := \{ j : \lim  | a_i - a_j | =0 \} . $$
Notice that, 
\begin{itemize}
\item $ i \in A ^ i$ and  $ A ^ i $ contains at least two elements. 
\item For $j \in A ^ i$, it holds that $ A ^ i = A ^ j $.
\item For each $j \in A ^ i$ and $k \notin A ^ i $, it holds: $ | a_j - a_k | \geq c  $ for some positive constant $c$.
\end{itemize}
Now, it is easy to see that 
$$ \{ 1, \cdots , N\} = \bigcup A ^ i  := A_1 \cup \cdots \cup A_q $$
with disjoint sets $ A_k $. 
Notice that 
$$ \e_{jk} \leq   \frac{c }{ \l^{n-2}} \qquad \mbox{ and } \qquad \frac{1}{\l_j } \Big| \frac{ \partial \e_{jk}}{ \partial a_j } \Big| \leq \frac{ c }{ \l ^{n-1} } \qquad \forall \, \, j \in A_i , \, \, \forall \, \, k \in A_\ell \mbox{ with } \ell \neq i . $$
Hence we are able to rewrite the equation \eqref{EAi1}, for $ i \in A _ \ell $ without the indices $k \notin A_ \ell $ by taking all these indices in the remainder term. Thus, we can repeat the argument done in Subsection 6.2 by taking $ A_ \ell $ instead of the set $ \{1, \cdots, N \}$ . 
Hence the proof is completed
\end{description}

\end{document}